\documentclass[preprint,12pt]{elsarticle}
\usepackage{mathtools}
\usepackage{amssymb}
\usepackage{amsthm}
\usepackage{ytableau}
\usepackage{eucal}
\usepackage{bm}
\usepackage{subcaption}

\newcommand{\N}{\mathbb{N}}
\newcommand{\la}{\lambda}
\renewcommand{\H}{\mathcal{H}}
\newcommand{\D}{\mathbf{D}}
\newcommand{\EMD}{\mathbf{EMD}}
\newcommand{\x}{\mathbf{x}}
\newcommand{\y}{\mathbf{y}}
\newcommand{\fraksl}{\mathfrak{sl}}
\newcommand{\Y}{{\rm Y}}
\newcommand{\scrY}{\mathcal{Y}}
\newcommand{\qbin}[2]{\genfrac[]{0pt}{}{#1}{#2}_q}
\newcommand{\qbininv}[2]{\genfrac[]{0pt}{}{#1}{#2}_{q^{-1}}}
\newcommand{\bla}{\bm{\la}}
\DeclareMathOperator{\PP}{PP}

\newtheorem{theorem}{Theorem}[section]
\newtheorem{corollary}[theorem]{Corollary}
\newtheorem{lemma}[theorem]{Lemma}
\newtheorem{proposition}[theorem]{Proposition}

\theoremstyle{definition}
\newtheorem{definition}[theorem]{Definition}
\newtheorem*{remark}{Remark}
\newtheorem{example}[theorem]{Example}
\newtheorem{problem}[theorem]{Problem}

\journal{\phantom{m}}

\begin{document}

\begin{frontmatter}



\title{Comparing weighted difference and earth mover's distance via Young diagrams}


\author[inst1]{William Q. Erickson}

\affiliation[inst1]{organization={Baylor University},
            addressline={1410 S.~4th St.}, 
            city={Waco},
            postcode={76706}, 
            state={TX},
            country={USA}}

\begin{abstract}
  We consider two natural statistics on pairs of histograms, in which the $n$ bins have weights $0, \ldots, n-1$. The difference ($\D$) between the weighted totals of the histograms is, in a sense, refined by the earth mover's distance ($\EMD$), which measures the amount of work required to equalize the histograms.  We were recently surprised, however, by how little the $\EMD$ actually does refine $\D$ in certain real-world applications, which led to the main problem in this paper: what is the probability that $\EMD = |\D|$?  We derive a formula for this probability, as well as the expected value of $|\D|$, via the combinatorics of Young diagrams and plane partitions.  We then generalize our results to an arbitrary number of histograms, where we realize this higher-dimensional $\D$ as distance on the Type-A root lattice.
\end{abstract}


\end{frontmatter}


\section{Introduction}

In this paper, we use the combinatorics of Young diagrams and plane partitions to resolve a statistical curiosity that arose from our analysis of students' grade distributions.  In an effort to assess the semesters during and after the period of virtual learning starting in 2020, the author joined a small group of department colleagues in studying grade distributions from the last several years of mathematics courses. In comparing two histograms of grades, the obvious statistic to track was the difference in GPA's (i.e., the weighted averages).  When tracking a single group of students across time, we could dispense with the average, and simply look at the \emph{weighted difference}: the signed difference between the total grade points earned in each class (where, for example, a grade of A equals 4 points, a B equals 3 points, etc.).  We will write $\D$ to denote this weighted difference.  Intuitively, $\D$ measures the minimum number of ``moves'' required to equalize the two weighted totals, where a ``move'' means moving one student up or down by one letter grade.

A statistician acquaintance then suggested that we study not only the weighted difference $\D$, but also the \emph{earth mover's distance} ($\EMD$).  Intuitively, the EMD measures the minimum number of ``moves'' required to equalize the two histograms themselves (not just their weighted totals), and is a true metric on the space of histograms.  Obviously it requires at least as many moves to equalize the histograms as it does merely to equalize their weighted totals; hence $\EMD$ is a sort of ``refinement'' of $\D$, in the sense that $\EMD$ always captures $\D$ along with some additional difference.  In general, then, one can certainly expect more information, and a different type of information, from studying $\EMD$ in tandem with $\D$.

\begin{figure}[ht]
    \centering
    \includegraphics[width=.7\linewidth]{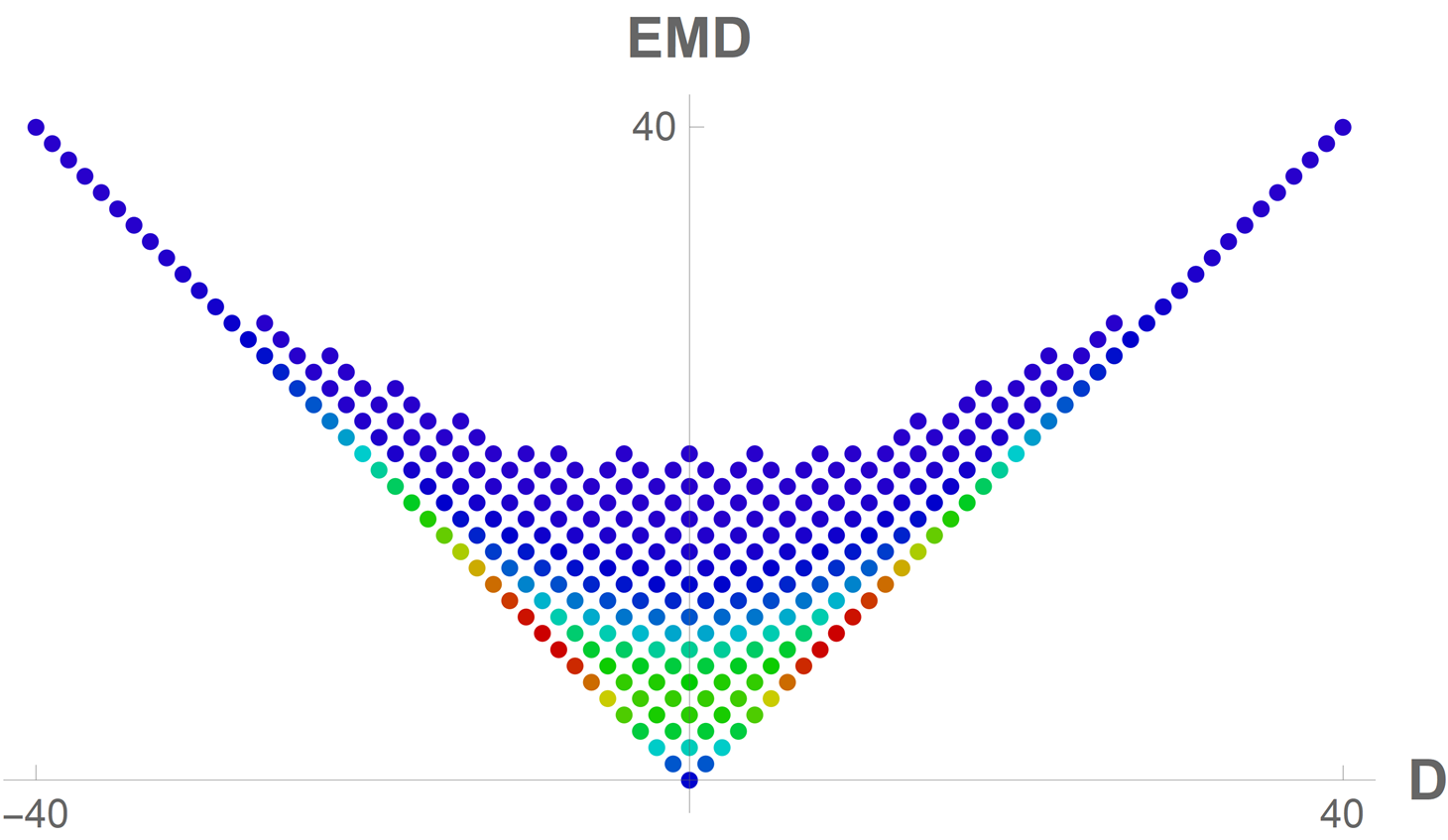}
    \caption{Example of a frequency plot for $\D$ vs.\ $\EMD$.  There are just over one million ordered pairs of histograms with $10$ data points and $5$ bins. The color of each dot corresponds to the number of pairs with given $\D$ and $\EMD$ (dark blue being the lowest, and red the highest).  Note the high concentration along the diagonals where $\EMD = |\D|$.}
    \label{fig:intro example}
\end{figure}

But after some initial excitement, we began to realize that in our data sets, $\EMD$ was almost always extremely close, or even equal, to the absolute value of $\D$. We were both disappointed and intrigued.  It is a fact that the EMD has proved tremendously powerful in an ever-widening range of applications, including computer vision \cite{rubner}, physics \cite{komiske}, cosmology \cite{frisch}, political science \cite{lupu}, epidemiology \cite{melnyk}, and many others.  Part of the issue in our case was the one-dimensional nature of the grading scale; another was the fact that the vast majority of grade distributions in real life approximate a normal distribution.  Still, even on our one-dimensional histograms, it is quite possible in theory for $\EMD$ and $|\D|$ to differ significantly (by roughly half the product of the number of data points and the number of bin dividers). Thus the following question emerged among our group: for a pair of histograms chosen uniformly at random, what is the probability that $\EMD = |\D|$?  (See Figure~\ref{fig:intro example}, which gives a sense of the concentration of histogram pairs with this property.)

The main result in this paper (Theorem~\ref{theorem:EMDequalsD}) is the answer to that question.  Our strategy is to view histograms as Young diagrams, in such a way that $\D$ is simply the size difference between two diagrams, while $\EMD$ is the symmetric difference. We then relate the two statistics by enumerating certain plane partitions.  We can quite easily extend this discrete result to the continuous setting of probability distributions (Corollary~\ref{cor:continuous}), where it turns out that the desired probability equals 2 divided by the number of histogram bins.  Our secondary result (Theorem~\ref{theorem:expd}) is a formula for the expected value of $|\D|$.  (The expected value of $\EMD$ was derived recursively in~\cite{bw}, and an explicit formula for probability distributions was subsequently obtained in~\cite{FV}.)  

In Section~\ref{section:general d}, we generalize the problem to an arbitrary number of histograms, rather than only two at a time.  Here the interplay between the Young diagrams and plane partitions is predictably more complicated, but nonetheless allows us  in Section~\ref{section:general d results} to write down generalizations of our two main results above.  We briefly conclude in Section~\ref{section:RepThy} by viewing our results from the perspective of representation theory, which raises a natural question for future research (Problem~\ref{problem:decompose}).

Certain passages in this paper are adapted from the author's thesis~\cite{EricksonThesis}.

\textbf{Acknowledgment:} We would like to thank the members of our ``2020 Hindsight'' group, especially Jeb Willenbring and Rebecca Bourn, for proposing the intensive grade analysis which ultimately motivated the problem in this paper.

\section{The weighted difference and the EMD}\label{sec:statistics}

Before shifting to a combinatorial perspective, we briefly describe the usual statistical approach to $\D$ and $\EMD$.  We consider  discrete histograms on a one-dimensional feature space, in which the bins are the elements of the set $[n]\coloneqq\{1,\ldots,n\}$.  Therefore we will use the term \emph{histogram} to refer to an element of $\N^n$ (with the convention $0 \in \N$). We will denote a histogram by the notation
\[
\la = (\la_1,\ldots,\la_n) \in \N^n,
\]
where $\la_i$ is the number of data points in bin $i$.  For positive integers $m \geq 1$ and $n \geq 2$, we define the set
\[
\H(m,n)\coloneqq \Big\{ \text{histograms with $m$ data points and $n$ bins}\Big\}.
\]
(These are also known as weak integer compositions of $m$ into $n$ parts.  We use the letter $m$ to evoke the ``mass'' from the original earth mover's problem \cite{monge}.)  It will be most natural to regard bin 1 as the ``best/highest'' and bin $n$ as the ``worst/lowest.''  Therefore we define the \emph{weighted total} of a histogram $\la$ as follows:
\begin{equation}\label{tdefinition}
t(\la) \coloneqq \sum_{i = 1}^{n} (n-i) \la_i.
\end{equation}
We define the \emph{weighted difference} $\D$ between two histograms $\la,\mu \in \H(m,n)$ to be
\begin{equation}
\label{D definition d=2}
    \D(\la,\mu) \coloneqq t(\la) - t(\mu),
\end{equation}
i.e., the signed difference between their weighted totals.  We observe that $\D$ can also be understood as the minimum amount of work required to equalize the two weighted totals, where one unit of \emph{work} is defined as moving one data point (in either histogram) up or down by exactly one bin.
\begin{example}
Suppose that $\la$ represents a grade distribution among $m$ students.  If $n=5$, and we regard the letter grades A, B, C, D, F as the bins $1, \ldots, 5$, then $t(\la)$ equals the total grade points earned, under the American system whereby A = 4, B = 3, C = 2, D = 1, and F = 0.  Given two grade distributions $\la,\mu \in \H(m,5)$, we are often interested in $\D(\la,\mu)/m$, the difference between the two GPA's.
\end{example}

Proceeding to the earth mover's distance (EMD), we consider $\la,\mu \in \H(m,n)$.  Intuitively, their EMD is the minimum amount of work (as defined above) required to transform $\la$ into $\mu$.  First we define the \emph{cost} (also known as ``ground distance'') between bins $i$ and $j$ to be their distance $|i-j|$.  We then have an $n \times n$ cost matrix $C$, where
$C_{ij} \coloneqq |i-j|$. For example, if $n=5$, then
\begin{equation*}
C = \begin{bmatrix}
    0 & 1 & 2 & 3 & 4\\
    1 & 0 & 1 & 2 & 3\\
    2 & 1 & 0 & 1 & 2\\
    3 & 2 & 1 & 0 & 1\\
    4 & 3 & 2 & 1 & 0
    \end{bmatrix}_{\textstyle .}
\end{equation*}
Note that $C_{ij}$ is the taxicab distance (i.e., the $\ell_1$-distance) from $(i,j)$ to the main diagonal of $C$.  (This is the property we will generalize for the higher-dimensional cost function in Section~\ref{section:general d}.)  

The EMD is defined in \cite{rubner} by the solution to a certain linear programming problem (traditionally named after Gaspard Monge \cite{monge}, or various combinations of Monge, Kantorovich, Hitchcock, and Koopman).  Specifically, we want to find a real $n \times n$ matrix $F$ in order to solve the following \emph{Monge problem}:
\begin{align}
    \text{Minimize} \quad & \sum_{\mathclap{i,j=1}}^n C_{ij} F_{ij}, \label{work}\\
    \text{subject to} \quad & F_{ij} \geq 0 & \text{for all } 1 \leq i,j \leq n, \label{nonneg}\\
    \text{and} \quad & \sum_{j=1}^n F_{ij} = \la_i & \text{for each }1\leq i \leq n, \label{rows}\\
    \text{and} \quad & \sum_{i=1}^n F_{ij} = \mu_j & \text{for each } 1 \leq j \leq n. \label{cols}
\end{align}
The matrix $F$ is called a \emph{flow matrix} or a \emph{transport plan} in transport theory literature.  (One can also view $F$ as a contingency table with $\la$ and $\mu$ as its margins.)  The idea is that any such $F$ encodes a way to transform $\la$ into $\mu$, by expending the amount of work in~\eqref{work}: namely, $F_{ij}$ equals the number of data points inside $\la$ that must ``flow'' from bin $i$ to bin $j$.  Hence \eqref{nonneg} any candidate $F$ must have nonnegative entries; \eqref{rows} the $i$th row sum (i.e., number of data points flowing out of bin $i$) must equal the number of data points in bin $i$ of $\la$; likewise, \eqref{cols} the $j$th column sum (the number of data points flowing into bin $j$) must equal the number of data points in bin $j$ of $\mu$.  Let $\mathcal F_{\la\mu}$ be the set of all flow matrices $F$ satisfying constraints \eqref{nonneg}--\eqref{cols}.

Once we have found an optimal solution $F$ to the Monge problem, the \emph{earth mover's distance} ($\EMD$) is defined to be the objective quantity in \eqref{work}; in other words,
\begin{equation}
    \label{emddef}
    \EMD(\la,\mu) \coloneqq \inf_{F \in \mathcal F_{\la\mu}} \left\{\sum_{i,j} C_{ij} F_{ij}\right\}.
\end{equation}
As a trivial example, consider the pair $(\la,\la)$.  Then the diagonal matrix $F = \operatorname{diag}(\la_1,\ldots,\la_n)$ lies in $\mathcal F_{\la\la}$, and we have $F_{ij} = 0$ whenever $i \neq j$.  Since $C_{ii}=0$, the total work in \eqref{work} equals $\sum_i C_{ii} \la_i = 0$, and so $\EMD(\la,\la)=0$, just as we would expect.

\begin{remark}
In many sources such as \cite{rubner}, the EMD is actually defined to be the normalization of our EMD, via division by $m$.  Here we omit normalizing so as to keep the combinatorics front and center.  Likewise, the EMD definition \eqref{emddef} certainly remains valid when the histogram components are allowed to be nonnegative \textit{real} numbers.  Hence, one could easily study probability distributions on $[n]$ rather than histograms.  Although we restrict our attention in this paper to discrete histograms, our approach nevertheless leads to a result (Corollary~\ref{cor:continuous}) for probability distributions, via normalizing by $1/m$ and then taking the limit as $m \rightarrow \infty$.

The EMD can also be defined to allow for \emph{partial matching}, i.e., histograms with different numbers of data points; for combinatorial approaches to partial matching, see the author's paper \cite{EricksonAStat} and thesis \cite{EricksonThesis}.  As explained in \cite{bickel}, the allowance for partial matching is precisely what distinguishes the EMD from the 1-Wasserstein (also known as Mallows) distance.
\end{remark}

By definition, computing the EMD depends upon first solving the Monge problem \eqref{nonneg}--\eqref{cols}, the difficulty of which varies with the choice of cost matrix $C$.  In our case, where we have $C_{ij} = |i-j|$, the Monge problem can be solved  via an $O(2n)$-time greedy algorithm known known as the ``northwest corner rule.''  (For details, see Hoffman \cite{hoffman}, who defined a certain \emph{Monge property} in cost matrices which is equivalent to a solution by the northwest corner rule.  It is not hard to see that our particular $C$ in this paper has this Monge property; see \cite[Lemma 2.4]{EricksonThesis}.)  

Recently, Bourn and Willenbring \cite{bw} improved upon the northwest corner rule by their use of the Robinson--Schensted--Knuth (RSK) correspondence, a ubiquitous combinatorial phenomenon with applications in many fields of mathematics. Translating the ideas in~\cite{bw} into the language of transport theory, we can summarize as follows.  First, write out the data points in $\la$ and $\mu$ in ascending order:
\begin{align}
\label{histograms to data points}
\begin{split}
        \la & \leadsto x_1 \ldots x_m = \underbrace{1 \ldots 1}_{\la_1} \underbrace{2 \ldots 2}_{\la_2} \ldots\ldots \underbrace{n \ldots n}_{\la_n},\\
        \mu & \leadsto y_1 \ldots y_m = \underbrace{1 \ldots 1}_{\mu_1} \underbrace{2 \ldots 2}_{\mu_2} \ldots\ldots \underbrace{n \ldots n}_{\mu_n}.
        \end{split}
    \end{align}
Then the RSK-image of these two sequences is precisely the optimal flow matrix $F$ for $\la$ and $\mu$. It follows from the RSK construction~\cite[p.~718]{knuth} that the nonzero entries of $F$ can be read off directly from the columns of the two-row array $\left(\begin{smallmatrix} x_1 & \ldots & x_m\\y_1 & \ldots & y_m\end{smallmatrix}\right)$.  Hence we can compute the EMD directly from the pairs $(x_j, y_j)$, thus bypassing the matrix $F$ entirely:
\begin{equation}
    \label{EMD compute RSK d=2}
    \EMD(\la,\mu) = \sum_{j=1}^m C_{x_j, y_j} = \sum_{j=1}^m |x_j - y_j|.
\end{equation}

\begin{example}\label{exampleearly}
Let $\la = (1,7,0,1,1)$ and $\mu = (4,0,1,3,2)$.  Then $\la,\mu \in \H(10,5)$.  We first arrange their data points as in~\eqref{histograms to data points}:
\begin{align*}
    \la &\leadsto 1222222245\\
    \mu & \leadsto 1111344455
\end{align*}
Then we simply sum the costs of the 10 columns $\binom{x_j}{y_j}$ to obtain
\begin{align*}
\EMD(\la,\mu) &= C_{11} + 3( C_{21}) + C_{23} + 3(C_{24})  + C_{45} + C_{55} \\
&= 0 + 3(0) + 1 + 3(2) + 1 + 0 = 11.
\end{align*}
Although there was no need to refer to the optimal flow matrix $F$, we can easily see that our computation above is equivalent to finding $\sum_{i,j} C_{ij} F_{ij}$ where
\[
F = \left[\begin{smallmatrix}
1 & 0 & 0 & 0 & 0\\
3 & 0 & 1 & 3 & 0\\
0 & 0 & 0 & 0 & 0\\
0 & 0 & 0 & 0 & 1\\
0 & 0 & 0 & 0 & 1
\end{smallmatrix}\right]
\]
is the RSK-image of our two sequences of data points.  The reader can check that $F$ encodes the instructions to transform $\la$ into $\mu$ in 11 moves.
\end{example}

\section{$\D$ and $\EMD$ in terms of Young diagrams}

In this section we realize each histogram $\la \in \H(m,n)$ as a Young diagram $\Y(\la)$ fitting inside an $m \times (n-1)$ rectangle.  A \emph{Young diagram} is a finite arrangement of rows of boxes, justified along the left and top edges, whose row lengths are weakly decreasing from top to bottom.  The number of boxes is called the \emph{size} of the Young diagram, denoted by vertical bars $| \cdot |$.  We want the size $|\Y(\la)|$ to equal the weighted total $t(\la)$, which suggests the following construction.  Write out the data points of $\la$ as in~\eqref{histograms to data points}, then replace each data point $x_j$ with its contribution $n-x_j$ to the weighted total.  The resulting sequence gives the row lengths of $\Y(\la)$. 
 For example:
 \[
 \ytableausetup{smalltableaux,centertableaux}
 \la = (2,1,0,3,1) \leadsto 1124445 \leadsto 4431110 \leadsto \Y(\la) = \ydiagram{4,4,3,1,1,1}
 \]
 We observe that $t(\la) = |\Y(\la)| = 14$.  In general, then, for $\la,\mu \in \H(m,n)$, we have $\D(\la,\mu) = |\Y(\la)| - |\Y(\mu)|$.

If we regard each Young diagram as a subset of $(\mathbb{Z}_{>0})^2$, with $(1,1)$ corresponding to the upper-left box, then it makes sense to write $Y_1 \supseteq Y_2$ in the usual sense of containment.  It also makes sense to speak of the \emph{symmetric difference} of two Young diagrams, where as usual the symmetric difference $\triangle$ of two sets $X$ and $Y$ is
$$
X \triangle Y \coloneqq (X \cup Y) \backslash (X \cap Y),
$$
i.e., the set of elements contained in exactly one of $X$ and $Y$.  Note that when $X$ and $Y$ are finite,
\begin{equation}
\label{sym diff = |D| if contain}
    |X \triangle Y| = \Big||X| - |Y|\Big| \iff X \supseteq Y \text{ or } Y \supseteq X.
\end{equation} 
The following proposition is essentially a combinatorial formulation of the fact \cite{delon} that the EMD of two distributions equals the $\ell_1$-distance between their cumulative distributions.  (Note that the $i$th column of $\Y(\la)$ counts the number of data points in the last $n-i$ bins.)

\begin{proposition}\label{prop:EMDisSymmDiff}
Let $\la,\mu \in \H(m,n)$.  Then $\EMD(\la,\mu) = \big|\Y(\la) \triangle \Y(\mu)\big|$.
\end{proposition}

\begin{proof}
To obtain the symmetric difference of the two diagrams, we can take the disjoint union over each of the rows.  By our construction of the Young diagrams, the $j$th row of $\Y(\la)$ has length $n-x_j$, and the $j$th row of $\Y(\mu)$ has length $n-y_j$, where $x_j$ and $y_j$ are the data points in~\eqref{histograms to data points}.  Therefore
\[
\big|\Y(\la) \triangle \Y(\mu)\big| = \sum_{j=1}^m |(n-x_j) - (n-y_j)| = \sum_j |x_j - y_j|,
\]
which by~\eqref{EMD compute RSK d=2} equals $\EMD(\la,\mu)$.
\end{proof}

\begin{example}
We return to the histograms from Example \ref{exampleearly}, namely $\la = (1,7,0,1,1)$ and $\mu = (4,0,1,3,2)$.  To obtain the row lengths of the Young diagrams, we have
\begin{align*}
\la &\leadsto 1222222245 \leadsto 4333333310\\
\mu &\leadsto 1111344455 \leadsto 4444211100
\end{align*}
and therefore
\[
\Y(\la) =\: \ytableausetup{centertableaux, smalltableaux} \ydiagram[*(red!30!white)]{4,3,3,3,3,3,3,3,1}
 \qquad \text{and} \qquad \Y(\mu) = \: \ydiagram[*(blue!30!white)]{4,4,4,4,2,1,1,1,0}.
 \]
 Incidentally, we see that $\D(\la,\mu)=26-21=5$.  We have shaded the diagrams red and blue to help visualize their symmetric difference, which we represent by the set of shaded boxes below:
 \[
\Y(\la) \triangle \Y(\mu) =  \begin{ytableau}
\phantom{1}&\phantom{1} & \phantom{1}&\phantom{1} \\
\phantom{1}&\phantom{1} & \phantom{1}& *(blue!30!white) \phantom{1}\\
\phantom{1}&\phantom{1} & \phantom{1}& *(blue!30!white) \phantom{1}\\
\phantom{1}&\phantom{1} & \phantom{1}& *(blue!30!white) \phantom{1}\\
\phantom{1} & \phantom{1}& *(red!30!white) \phantom{1}\\
\phantom{1} & *(red!30!white)\phantom{1}& *(red!30!white) \phantom{1}\\
\phantom{1} & *(red!30!white)\phantom{1}& *(red!30!white) \phantom{1}\\
\phantom{1} & *(red!30!white)\phantom{1}& *(red!30!white) \phantom{1}\\
*(red!30!white)\phantom{1}\\
\end{ytableau}
 \]
 Counting the shaded boxes, we find that $|\Y(\la) \triangle \Y(\mu)| = 11$, which agrees with our result for $\EMD(\la,\mu)$ from Example \ref{exampleearly}.
\end{example}

Let $\scrY(a,b)$ denote the set of all Young diagrams fitting inside an $a \times b$ rectangle.  Clearly $\la \in \H(m,n)$ if and only if $\Y(\la) \in \scrY(m,\:n-1)$.  The cardinality of each of these sets is well known to be $\binom{m+n-1}{m}$, which we can see at once by viewing a histogram as $m$ indistinguishable data points along with $n-1$ dividers between bins, since choosing $m$ locations for the data points uniquely determines a histogram.  In fact, by transposing Young diagrams, we have the following equalities as well:
\begin{align*}
    \binom{m+n-1}{m} = \# \H(m,n) &= \# \scrY(m, \: n-1) \\
    &= \# \scrY(n-1,\:m) = \# \H(n-1,\:m+1).
\end{align*}
Since the size and symmetric difference of Young diagrams are unaffected by transposing the diagrams, it follows that any statistical results involving $\D$ and $\EMD$ on $\H(m,n)$ also hold true on $\H(n-1, m+1)$.  We might give this duality the name ``Hermite reciprocity for histograms,'' since at a deeper level, the relevant identity is actually the  equivalence $S^m(\mathbb{C}^n) \cong S^{n-1}(\mathbb{C}^{m+1})$ of representations for the Lie algebra $\fraksl(2,\mathbb{C})$.  (The $\D$-generating function $f(q)$ in Theorem~\ref{theorem:expd} below can be viewed as a character for $\fraksl(2,\mathbb{C})$; see Section~\ref{section:RepThy} of this paper, and~\cite[p.~110]{EricksonThesis}.)

\section{Plane partitions}
\label{section:plane partitions}

We will prove our main results using a generalization of Young diagrams known as plane partitions.  A \emph{plane partition} is a Young diagram in which the boxes are filled with positive integers, weakly decreasing along each row and column.  A plane partition is visualized three-dimensionally by stacking cubes on a Young diagram such that the height above each box is its integer entry; for this reason, a plane partition whose underlying Young diagram has at most $x$ rows and $y$ columns, with entries at most $z$, is said to fit inside an $x \times y \times z$ box.  The number of such plane partitions is
\begin{equation}
\label{planepartbox}
\PP(x,y,z)= \prod_{i=1}^x \prod_{j=1}^y \prod_{k=1}^z \frac{i+j+k-1}{i+j+k-2.}
\end{equation}
This formula was proved by MacMahon in \cite{MacMahon}; see also \cite[\S 7.20--22]{stanley2} for an expansive treatment of the combinatorics of plane partitions. 

A plane partition fitting inside an $m \times (n-1) \times d$ box can be viewed as a chain of Young diagrams $Y_1, \ldots, Y_d \in \scrY(m, n-1)$ with respect to containment.  Specifically, given a plane partition, $Y_i$ is the underlying diagram containing all entries greater than or equal to $i$. Conversely, given $Y_1 \supseteq \cdots \supseteq Y_d$, fill each box of $Y_1$ with the number of diagrams $Y_i$ containing that box.  Hence we have a bijective correspondence:
\begin{equation}
\label{bijection PP's}
    \left\{\quad\parbox{2.5cm}{\centering chains\\
    $Y_1 \supseteq \cdots \supseteq Y_d$\\ in $\scrY(m, n-1)$}\quad\right\} \longleftrightarrow \left\{\quad\parbox{3.5cm}{\centering plane partitions fitting inside an $m \times (n-1) \times d$ box}\quad\right\}
\end{equation}
The three-dimensional visualization of the plane partition is obtained by stacking the diagrams on top of each other, with $Y_1$ on bottom and $Y_d$ on top.

In the following lemma, which records three specializations of the formula~\eqref{planepartbox} for later use, we write the Pochhammer symbol $(x)_y \coloneqq \frac{\Gamma(x+y)}{\Gamma(x)} = x(x+1) \cdots (x+y-1)$ to denote the rising factorial. 

\begin{lemma}
\label{lemma:planepart}
We have
\begin{align}\setlength{\linewidth}{3cm}
    \PP(m, \: n-1, \: 1) &= \frac{(n)_m}{m!}, \label{PP1}\\[2ex]
    \PP(m, \: n-1, \: 2) &= \frac{(n)_m (n+1)_m}{m!(m+1)!}, \label{PP2}\\[2ex]
    \PP(m, \: n-1, \: 3) &= \frac{2(n)_m (n+1)_m (n+2)_m}{m! (m+1)! (m+2)!}. \label{PP3}
\end{align}
\end{lemma}

\begin{proof}

For \eqref{PP1}, rather than deriving directly from~\eqref{planepartbox}, we merely observe that $\PP(m, \: n-1, \: 1) = \#\scrY(m, n-1) = \binom{m+n-1}{m} = \frac{(n)_m}{m!}$.  For~\eqref{PP2}, we separate \eqref{planepartbox} into two main factors (one for $k=1$ and the other for $k=2$) and cancel the repeated factors:
\[
    \PP(x,y,2) = \left(\prod_{i,j} \frac{i+j}{i+j-1} \right) \left(\prod_{i,j} \frac{i+j+1}{i+j} \right)\\
    = \prod_{i,j} \frac{i+j+1}{i+j-1.}
\]
Now we use induction on $x = m$ and $y = n-1$ to prove that this product is equivalent to~\eqref{PP2}.  In the base case $x=y=1$, the product above is just the single factor $\frac{1+1+1}{1+1-1}=3$, and likewise the expression in~\eqref{PP2} is $\frac{2 \cdot 3}{1!2!}=3$.  
To induct on $x$, assume that~\eqref{PP2} holds for $x-1 = m-1$ and $y = n-1$.  Then we have
\begin{align*}
    \PP(m, \: n-1, \: 2) &= \prod_{i=1}^{m} \prod_{j=1}^{n-1} \frac{i+j+1}{i+j-1}\\
    &= \prod_{j=1}^{n-1} \frac{m+j+1}{m+j-1} \cdot \PP(m-1,n-1,2)\\[2ex]
    &= \frac{(m+n-1)(m+n)}{m(m+1)} \cdot \frac{(n)_{m-1} (n+1)_{m-1}}{m!(m-1)!}\\
    &= \frac{(n)_m (n+1)_m}{m!(m+1)!}
\end{align*}
as in~\eqref{PP2}.  The induction on $y = n-1$ is identical.  The proof for~\eqref{PP3} is carried out in the same way, so we omit the intermediate simplifications.
\end{proof}

\section{Main results for histogram pairs}

In this section, we use the combinatorics of Young diagrams and plane partitions to answer the two statistical questions motivating the paper: assuming the uniform distribution on $\H(m,n) \times \H(m,n)$,
\begin{enumerate}
    \item With what probability does $\EMD = |\D|$?
    \item What is the expected value of $|\D|$?
\end{enumerate}



\begin{theorem}\label{theorem:EMDequalsD}
\label{emc=d}
Let $(\la,\mu)\in \H(m,n) \times \H(m,n)$ be chosen uniformly at random.  Then 
\begin{equation}
    \label{formula EMD=D for d=2}
    \mathbb P\Big(\EMD(\la,\mu) = |\D(\la,\mu)|\Big) =   \frac{2(m+n)}{n(m+1)} - \frac{m!}{(n)_m}.
\end{equation}
\end{theorem}

\begin{proof} 
By~\eqref{sym diff = |D| if contain}, we have that $\EMD(\la,\mu)=|\D(\la,\mu)|$ if and only if one of the two diagrams $\Y(\la)$ and $\Y(\mu)$ contains the other. By the bijection in~\eqref{bijection PP's}, the number of such \emph{un}ordered pairs equals $\PP(m, \: n-1, \: 2)$.  To find the number of such \emph{ordered} pairs, we multiply by 2 and then correct for the pairs $(\la,\la)$ by subtracting $\#\H(m,n) =\binom{m+n-1}{m} = \frac{(n)_{m}}{m!}$. Hence the number of ordered pairs for which $\EMD = |\D|$ is
\[
2\cdot \PP(m, \: n -1, \:2) - \frac{(n)_{m}}{m!}.
\]
We have a formula from \eqref{PP2} in Lemma~\ref{lemma:planepart} for the first term. 
To obtain the desired probability, we divide by $\#\H(m,n)^2$, and after simplifying we obtain
\begin{align*}
    \quad \Bigg(2\cdot \overbrace{\frac{(n)_m (n+1)_m}{m!(m+1)!}}^{\PP(m,\: n-1, \: 2)} - \frac{(n)_{m}}{m!}\Bigg) \Bigg/ \bigg(\!\!\overbrace{\frac{(n)_{m}}{m!}}^{\#\H(m,n)}\!\!\bigg)^2 &= \frac{2m!(n+1)_m}{(n)_m (m+1)!} - \frac{m!}{(n)_m}\\[2ex]
    &= \frac{2(m+n)}{n(m+1)} - \frac{m!}{(n)_m}.
    \qedhere
\end{align*}
\end{proof}

 See Table \ref{table:EMDequalsD}, which records this probability for selected values of $m$ and $n$.  This result also translates nicely into the setting of probability distributions, as follows. Consider the normalized versions of our two statistics:
\[
\widehat{\EMD}(\la,\mu)\coloneqq \frac{\EMD(\la,\mu)}{m} \quad \text{and} \quad \widehat{\D}(\la,\mu) \coloneqq \frac{\D(\la,\mu)}{m}.
\]
Then both statistics assume real values in the interval $[0,\: n-1]$. Both $\widehat{\EMD}$ and  $\widehat{\D}$ have natural analogues for the set $\mathcal{P}(n)$ of probability distributions on $[n]$, obtained by letting $m \rightarrow \infty$.  Hence by taking the limit of~\eqref{formula EMD=D for d=2} as $m \rightarrow \infty$, we arrive at a result for probability distributions, rather than for discrete histograms.  (The reader can compare the corollary below to the bottom row of Table~\ref{table:EMDequalsD}.)

\begin{corollary}
\label{cor:continuous}
Let $(X,Y) \in \mathcal{P}(n) \times \mathcal{P}(n)$ be chosen uniformly at random.  Then
\[
\mathbb{P}\Big(\widehat{\EMD}(X,Y) = |\widehat{\D}(X,Y)|\Big) = \frac{2}{n}.
\]
\end{corollary}

\begin{table}[ht]
    \centering
   
\begin{tabular}{|c|c|c|c|c|c|c|c|c|c|}
\hline
$m \backslash n$ & 2 & 3 & 4 & 5 & 6 & 7 & 8 & 9 & 10\\ \hline
 10&1&0.773 & 0.633 & 0.544 & 0.485 & 0.441 & 0.409 & 0.384 & 0.364 \\
 20&1&0.726 & 0.571 & 0.476 & 0.413 & 0.367 & 0.333 & 0.307 & 0.286 \\
 30&1&0.708 & 0.548 & 0.452 & 0.387 & 0.341 & 0.306 & 0.280 & 0.258 \\
40&1& 0.698 & 0.537 & 0.439 & 0.374 & 0.328 & 0.293 & 0.266 & 0.244 \\
50&1& 0.692 & 0.529 & 0.431 & 0.366 & 0.319 & 0.284 & 0.257 & 0.235 \\
60&1& 0.688 & 0.525 & 0.426 & 0.361 & 0.314 & 0.279 & 0.251 & 0.230 \\
 70&1&0.685 & 0.521 & 0.423 & 0.357 & 0.310 & 0.275 & 0.247 & 0.225 \\
 80&1&0.683 & 0.519 & 0.420 & 0.354 & 0.307 & 0.272 & 0.244 & 0.222 \\
 90&1&0.681 & 0.516 & 0.418 & 0.352 & 0.305 & 0.269 & 0.242 & 0.220 \\
 100&\phantom{10}1\phantom{10}&0.680 & 0.515 & 0.416 & 0.350 & 0.303 & 0.267 & 0.240 & 0.218 \\ 
 10,000 & 1 & 0.667 & 0.500 & 0.400 & 0.333 & 0.286 & 0.250 & 0.222 & 0.200 \\
\hline
\end{tabular}
    \caption[Values of $\mathbb P(\EMD = |\D|)$ on $\H(m,n) \times \H(m,n)$]{Values of $\mathbb P(\EMD = |\D|)$ on $\H(m,n) \times \H(m,n)$.  Obtained in Mathematica using Theorem \ref{theorem:EMDequalsD}.}
    \label{table:EMDequalsD}
\end{table}

Our second main question was to find the expected value of $|\D(\la,\mu)|$.  In order to express this result, we recall the $q$-analogue of the binomial coefficients $\binom{a}{b}$.  First we define $[a]_q \coloneqq 1 + q + \cdots + q^{a-1}$.  Then writing $[a]_q! \coloneqq [a]_q [a-1]_q \cdots [1]_q$, we define the \emph{$q$-binomial coefficient} 
\[
\qbin{a}{b} \coloneqq \frac{[a]_q!}{[b]_q![a-b]_q!},
\]
which is a polynomial in $q$ encoding combinatorial information central to this paper:
\begin{equation}\label{qbinyoungfact}
\qbin{m+n-1}{m} = \sum_{Y \in \scrY(m,\:n-1)} q^{|Y|} = \sum_{\la \in \H(m,n)} q^{t(\la)},
\end{equation} where the first equality is Proposition 1.7.3 in \cite{stanley1}.

We will borrow the \emph{plussing} operator $[\quad ]_+$ from the theory of economic modeling; this operator, also called a ``linear annihilation operator'' in \cite{hansen}, acts on a Laurent series in $q$ by annihilating negative powers of $q$:
$$\left[ \sum_{k=-\infty}^{\infty} c_k q^k\right]_+ \coloneqq \sum_{k=0}^\infty c_k q^k.
$$
See \cite[p.~14]{hansen}, and the original Wiener--Kolmogorov prediction formula dating back to \cite{wiener, kolmogorov}.

\begin{theorem}
\label{theorem:expd}
Let $(\la,\mu) \in \H(m,n) \times \H(m,n)$ be chosen uniformly at random.  Then we have
\[
\mathbb{E}\Big(\big|\D(\la,\mu)\big|\Big) = \frac{2\, g'(1)}{\binom{m+n-1}{m}_,^2}
\]
where $g(q) = [f(q)]_+$ and $f(q) = \qbin{m+n-1}{m} \qbininv{m+n-1}{m}$.    
\end{theorem}
\begin{proof}
By \eqref{qbinyoungfact}, we have
\begin{align*}
    f(q) \coloneqq \qbin{m+n-1}{m} \qbininv{m+n-1}{m}&= \sum_{\la \in \H(m,n)}q^{t(\la)} \sum_{\mu \in \H(m,n)}q^{-t(\mu)}\\
    &= \sum_{\la,\mu} q^{t(\la)-t(\mu)}\\
    &= \sum_{\la,\mu} q^{\D(\la,\mu)}.
\end{align*}
Now applying the plussing operator, we have
\[
    g(q) \coloneqq [f(q)]_+ = \sum_{t(\la)\geq t(\mu)}q^{\D(\la,\mu)} = \sum_{k = 0}^{m(n-1)} c_k\:q^k
\]
where $c_k = \#\{(\la,\mu) \mid \D(\la,\mu) = k\}$.  Hence
\begin{align*}
    g'(q) &= \sum_{k = 0}^{m(n-1)} (c_k\cdot k)q^{k-1}
\end{align*}
and therefore, since by the symmetry $(\la,\mu) \leftrightarrow (\mu,\la)$ we must have $c_k = c_{-k}$, we obtain
\begin{align*}
    g'(1) &= \sum_{k = 0}^{m(n-1)} c_k \cdot k\\ 
    &= \sum_{k = 0}^{m(n-1)} c_{-k} \cdot k\\
    &= \frac{1}{2} \cdot \sum_{k = 1}^{m(n-1)} \#\Big\{(\la,\mu) \: \Big\vert \: |\D(\la,\mu)| = k \Big\}\cdot  k.
\end{align*}
Upon multiplying both sides by 2, the right-hand side is just the expected value of $|\D(\la,\mu)|$ multiplied by $\#\H(m,n)^2$, so we are done.
\end{proof}
\begin{remark}
Similarly as in Corollary~\ref{cor:continuous}, we can normalize this result by substituting $\sqrt[m]{q}$ for $q$.  Then the asymptotic as $m \rightarrow \infty$ equals the expected value of $|\widehat \D|$ on $\mathcal{P}(n) \times \mathcal{P}(n)$.
\end{remark}

Using the combinatorial formula in Theorem \ref{theorem:expd}, we can quickly calculate expected values of $|\D|$ in Mathematica.  See Table \ref{table:expvald}, in which we have \emph{unit} normalized to obtain values between 0 and 1.  For example, from Table \ref{table:expvald}, we see that the expected GPA difference for two classes of 30 students equals $4(0.177)=0.708$.  We have double-checked the values in both Tables \ref{table:EMDequalsD} and \ref{table:expvald} by brute force in Mathematica, for sufficiently small $m$ and $n$.  Of course, brute force quickly becomes impractical, since (for example) $\H(100,10) \times \H(100,10)$ contains approximately $1.8 \times 10^{25}$ pairs; and yet the combinatorial approach allows us to compute both statistics over all these pairs in just a few seconds.

\begin{table}[]
    \centering

\begin{tabular}{|c|c|c|c|c|c|c|c|c|c|}
\hline
$m \backslash n$ & 2 & 3 & 4 & 5 & 6 & 7 & 8 & 9 & 10\\ \hline
 10 & 0.364 & 0.266 & 0.224 & 0.201 & 0.185 & 0.174 & 0.166 & 0.159 & 0.154 \\
 20 & 0.349 & 0.250 & 0.208 & 0.183 & 0.167 & 0.155 & 0.146 & 0.139 & 0.133 \\
 30 & 0.344 & 0.245 & 0.202 & 0.177 & 0.160 & 0.148 & 0.139 & 0.132 & 0.126 \\
 40 & 0.341 & 0.242 & 0.199 & 0.174 & 0.157 & 0.145 & 0.135 & 0.128 & 0.122 \\
 50 & 0.340 & 0.240 & 0.197 & 0.172 & 0.155 & 0.143 & 0.133 & 0.126 & 0.119 \\
 60 & 0.339 & 0.239 & 0.196 & 0.171 & 0.154 & 0.141 & 0.132 & 0.124 & 0.118 \\
 70 & 0.338 & 0.238 & 0.195 & 0.170 & 0.153 & 0.140 & 0.130 & 0.123 & 0.116 \\
 80 & 0.337 & 0.238 & 0.194 & 0.169 & 0.152 & 0.139 & 0.130 & 0.122 & 0.115 \\
90 & 0.337 & 0.237 & 0.194 & 0.168 & 0.151 & 0.139 & 0.129 & 0.121 & 0.115 \\
 100 & 0.337 & 0.237 & 0.193 & 0.168 & 0.151 & 0.138 & 0.128 & 0.121 & 0.114 \\ \hline
\end{tabular}

    \caption[Expected value of $|\D|$ on $\H(m,n) \times \H(m,n)$]{Expected value of $|\D|$ on $\H(m,n) \times \H(m,n)$, unit normalized by dividing by $m(n-1)$.  Obtained in Mathematica using the combinatorial formula in Theorem \ref{theorem:expd}.}
    \label{table:expvald}
\end{table}

\section{A generating function for $\D$ and $\EMD$}\label{section:genfunEMDandD}

We set aside this brief section to point out a known generating function for $\D$ and $\EMD$, which we reinterpret in terms of Young diagrams.  In \cite{kretschmann}, a recursive definition is given for a generating function that refines the EMD generating function proved in \cite{bw}.  Because of its recursive nature, we must formally consider histogram pairs with unequal numbers of bins ($n_1$ and $n_2$).  The EMD is still defined on $\H(m,n_1) \times \H(m,n_2)$ since we can append zeros so that all histograms are tuples of length $\max\{n_1,\:n_2\}$.  We do not consider the statistical meaning of the EMD in such a scenario; the unequal bin numbers are just a means to obtain a recursion.

\begin{proposition}[\cite{kretschmann}]
\label{kretschmanngenfun}
Formally, define 
\begin{equation}
    \label{H_nm}
H_{n_1,n_2} = H_{n_1,n_2}(q,x,y,z) \coloneqq \sum_{m=0}^\infty \left(\sum_{\substack{\la \in \H(m,n_1),\\ \mu \in \H(m,n_2)}}q^{\EMD(\la,\mu)}x^{t(\la)} y^{t(\mu)} \right)z^m.
\end{equation}
Then we have the recursion
\begin{equation} 
\label{H_nm recur}
H_{n_1,n_2}=\frac{H_{n_1-1,n_2}+H_{n_1,n_2-1}-H_{n_1-1,n_2-1}}{1-q^{C_{n_1,n_2}}x^{n_1-1}y^{n_2-1}z}
\end{equation}
where $H_{1,1}=\frac{1}{1-z}$, and $H_{0,n_2}=H_{n_1,0}=0$.
\end{proposition}
 Note that upon substituting $y = x^{-1}$, the indeterminate $x$ keeps track of $\D$.  Below we present a short alternate proof entirely in terms of Young diagrams.  Conversely, the change in perspective suggests that the generating function~\eqref{H_nm} may be of interest not only in statistics but also in combinatorics.
 
 \begin{proof}
 First we reinterpret the generating function \eqref{H_nm} in terms of Young diagrams, setting $X = \Y(\la)$ and $Y = \Y(\mu)$:
 \begin{equation*}
     H_{n_1,n_2}(q,x,y,z)=\sum_{m=0}^\infty \left(\sum_{\substack{X \in \scrY(m,\:n_1-1),\\ Y \in \scrY(m,\:n_2-1)}}q^{|X \triangle Y|}x^{|X|} y^{|Y|} \right)z^m.
 \end{equation*}
 For the base case $H_{1,1}$ it is clear that $\scrY(m,0)$ contains only the empty diagram, so the coefficient of $z^m$ has only one term; since the empty diagram has size 0, the exponents of $x$ and $y$ in this term are both 0.  Moreover, the symmetric difference of two empty diagrams is also empty, so the exponent of $q$ is also 0.  Hence the coefficient of $z^m$ in $H_{1,1}$ is $q^0 x^0 y^0=1$, so $H_{1,1}=1+z+z^2 +\cdots = \frac{1}{1-z}$.  
 For the second base case, in which $n_1$ or $n_2$ is 0, we have $\scrY(m,\:-1)=\varnothing$, since a Young diagram cannot have a negative number of columns.
 
 Now we consider $H_{n_1,n_2}$ in general. For any pair $(X,Y)$ in the inside sum, the diagram $X$ has a certain number $r_X$ of maximal rows (length $n_1-1$) at the top, and likewise $Y$ has a certain number $r_Y$ of maximal rows (length $n_2-1$) at the top.  Hence we can remove $r\coloneqq\min\{r_X,r_Y\}$ maximal rows from each diagram.  Since $z$ tracks the number of rows in each diagram, each removed row-pair had contributed 1 to the exponent of $z$.  Note that the size of the symmetric difference of a maximal row from $X$ and a maximal row from $Y$ is exactly $|n_1-n_2|=C_{n_1,n_2}$, so each of the removed row-pairs had contributed $C_{n_1,n_2}$ to the exponent of $q$.  Moreover, each removed row-pair had contributed the size $n_1-1$ to the exponent of $x$, and $n_2-1$ to the exponent of $y$.  Therefore we see that from each term in $H_{n_1,n_2}$ we can factor out $(q^{C_{n_1,n_2}}x^{n_1-1}y^{n_2-1}z)^r$ for some $r$, leaving us with a pair of diagrams \emph{at least one of which contains no maximal rows}; that is, at least one of them can fit inside a rectangle with one less column than before.  Hence, starting with $H_{n_1,n_2}$ and factoring out the reciprocal of the denominator in \eqref{H_nm recur}, the remaining terms are precisely those appearing in $H_{n_1-1,n_2}$ or in $H_{n_1,n_2-1}$.  After correcting for double-counting the terms appearing in both (namely $H_{n_1-1,n_2-1}$), we obtain the numerator of \eqref{H_nm recur}, and so we are done.
 \end{proof}
 
\section{Generalization of $\D$ and $\EMD$ to more than two histograms}
\label{section:general d}

In this section, we generalize $\D$ and $\EMD$ in order to compare an arbitrary number of histograms at a time, rather than only two.  We let $d$ denote the number of histograms under comparison.  Since we will work with $d$-tuples of histograms, we need a way to write down elements of the $d$-fold Cartesian product $\H(m,n)^d$, and so we use a boldface Greek letter 
$$\bla = (\la^1,\ldots,\la^d),
$$
where each $\la^i = (\la^i_1,\ldots,\la^i_n) \in \H(m,n)$.

\subsection{Generalization of $\D$}

As before with two histograms, we still want $\D$ to capture just the right amount of information, namely, each of the pairwise signed differences between weighted totals, but nothing more.  Therefore we propose
\[
\Lambda_d \coloneqq \mathbb{Z}^d / \mathbb{Z}(1, \ldots, 1)
\]
as the natural codomain of $\D$.  We write an element of $\Lambda_d$ using square brackets around any of its representatives in $\mathbb{Z}^d$.  For example, $[1,5,2] = [4,8,5] = [-1,3,0]$.

\begin{definition}
\label{def:Dgeneral}
Let $d\geq 2$.  We (re)define the \emph{weighted difference} to be the function $\D: \H(m,n)^d \longrightarrow \Lambda_d$ given by
\[
\D(\la^1,\ldots,\la^d) \coloneqq \big[t(\la^1),\ldots,t(\la^d)\big].
\]
\end{definition}

Note that each element of $\Lambda_d$ has a unique representative whose last coordinate is $0$.  This unique representative can then be truncated into a $(d-1)$-tuple.  With this convention in mind, our new definition of $\D$ agrees with the original definition \eqref{D definition d=2} in the $d=2$ case, since $\D(\la,\mu) = \big[t(\la), t(\mu)\big] = \big[ t(\la) - t(\mu),\:0\big] \leadsto t(\la) - t(\mu)$, just as before.

\begin{remark}
The reader may have recognized $\Lambda_d$ as the weight lattice (or the root lattice) in the root system of Type ${\rm A}_{d-1}$. 
 See Section~\ref{section:RepThy} for the interpretation of our results from the perspective of the representation theory of the Lie algebra $\fraksl(d,\mathbb{C})$.
\end{remark}

\begin{example}
Consider $m$ voters participating in a ranked-choice election between candidates $1, \ldots, d$.  Then for each candidate $i$, the results are represented by a histogram $\la^i \in \H(m,d)$, where $\la^i_j$ is the number of $j$th-place votes received by candidate $i$.  In this context, $t(\la^i)$ equals the \emph{Borda count} for candidate $i$.  In the Borda system, candidate $i$ wins the election if and only if the $i$th coordinate in any representative of $\D(\bla)$ is greater than all other coordinates.
\end{example}

\subsection{Generalization of $\EMD$}

Intuitively, we still want the EMD to count the minimum amount of work required to equalize all $d$ of the histograms.  We therefore adopt the same generalization as we did in~\cite{EricksonAStat}, which we summarize here.  The Monge problem~\eqref{work}--\eqref{cols} naturally generalizes to $d$ histograms, in which case $F$ is a $d$-dimensional array with side length $n$.  Instead of ``row'' and ``column'' sums in~\eqref{rows} and~\eqref{cols}, we simply require for $F$ that the sum of the entries in the $j$th hyperplane perpendicular to the $i$th coordinate axis must equal $\la^i_j$. 

Rather than a cost matrix, $C$ is now a $d$-dimensional cost array with side length $n$.  Just as in the $d=2$ case, we want the entry $C(\x)$ to be the taxicab distance from $\x$ to the main diagonal of $C$.  If $\x = (x_1, \ldots, x_d)$, then let $\widetilde{\x} = (\widetilde{x}_1, \ldots, \widetilde{x}_d)$ be the vector whose coordinates are the $x_i$ sorted in descending order.  It is shown in~\cite[Prop.~1]{EricksonAStat} that the taxicab distance to the diagonal can be computed as follows:  
\begin{equation}
    \label{C hard}
C(\x) = \sum_{i=1}^{\lfloor d/2 \rfloor} \widetilde{x}_i - \widetilde{x}_{d+1-i}.
\end{equation}
To reformulate this more conveniently, we define a special vector $\mathbf{u}$, which will be important for comparing $\D$ and $\EMD$:
\begin{equation}
\label{u}
\mathbf{u} \coloneqq \begin{cases}
    (\underbrace{1, \ldots, 1}_{d/2},\underbrace{-1, \ldots, -1}_{d/2}), & \text{$d$ even},\\[5ex]
    (\underbrace{1, \ldots, 1}_{(d-1)/2},0,\underbrace{-1, \ldots, -1}_{(d-1)/2}), & \text{$d$ odd}.
\end{cases}
\end{equation}
We now rewrite~\eqref{C hard} as 
\begin{equation}
    \label{cost}
    C(\x) = \mathbf{u} \cdot \widetilde{\x},
\end{equation}
where the dot denotes the standard dot product. For example, if $\x = (6,1,8,3,4)$, then we have
\[
C(\x) = \mathbf{u} \cdot \widetilde{\x} = (1,1,0,-1,-1) \cdot (8,6,4,3,1) = 8+6-(3+1) = 10.
\]
It turns out that the cost array $C$ retains the Monge property in any dimension $d$, and so the EMD can be computed via the obvious generalization of the RSK algorithm~\eqref{histograms to data points}--\eqref{EMD compute RSK d=2}.  Namely, arrange the data points from each $\la^i$ in ascending order
\begin{equation}
\label{histogram to data points general d}
\la^i \leadsto x^i_1 \ldots x^i_m = \underbrace{1 \ldots 1}_{\la^i_1} \underbrace{2 \ldots 2}_{\la^i_2} \ldots\ldots \underbrace{n \ldots n}_{\la^i_n}
\end{equation}
and then (bypassing the optimal flow array $F$) we simply sum the costs of the $m$ resulting $d$-tuples:
\begin{equation}
\label{EMD RSK general d}
\EMD(\bla) = \sum_{j=1}^m C(x^1_j, \ldots, x^d_j).
\end{equation}

\begin{example}
\label{example:RSK d general}
Set $d=4$, $n=5$, and $m=6$.  We will compute $\EMD(\bla)$, where $\bla = (\la^1, \ldots, \la^4)$ for the four histograms below:
\begin{align*}
    \la^1 &= (4,1,1,0,0) \leadsto  111123\\
    \la^2 &= (3,0,0,0,3) \leadsto  111555\\
    \la^3 &= (0,4,2,0,0) \leadsto  222233\\
    \la^4 &= (1,1,2,1,1) \leadsto  123345
\end{align*}
Using~\eqref{cost} to compute the costs of each $(x^1_j, \ldots, x^4_j)$, we have
\begin{align*}
    C(1,1,2,1) &= \mathbf{u} \cdot (2,1,1,1) = 1,\\
    C(1,1,2,2) &= \mathbf{u} \cdot (2,2,1,1) = 2,\\
    C(1,1,2,3) &= \mathbf{u} \cdot (3,2,1,1) = 3,\\
    C(1,5,2,3) &= \mathbf{u} \cdot (5,3,2,1) = 5,\\
    C(2,5,3,4) &= \mathbf{u} \cdot (5,4,3,2) = 4,\\
    C(3,5,3,5) &= \mathbf{u} \cdot (5,5,3,3) = 4.
\end{align*}
Summing these costs, we conclude $\EMD(\bla) = 1+2+3+5+4+4 = 19$.
\end{example}

\subsection{Generalized $\EMD$ via Young diagrams}

Typically in the literature, the definition of symmetric difference $\triangle$ is extended associatively to any finite number of sets; as a result, the symmetric difference contains those elements that are contained in an odd number of the original sets.  In order to translate the EMD into the language of Young diagrams, however, we introduce a different generalization in terms of multisets.  We write $x^a$ to denote a multiset element $x$ with multiplicity $a$.  Note that $\{x^0\} =\varnothing$.

\begin{definition}\label{def:usd}
The generalized symmetric difference of finite sets $X_1,\ldots,X_d$ is the multiset
\[
\blacktriangle(X_1,\ldots,X_d)\coloneqq \bigcup_{k=0}^d \{ x^{\min\{k,d-k\}} \mid \text{$x$ is contained in exactly $k$ of the $X_i$}\}.
\]

\end{definition}

It follows from the definition that if the $X_i$ are all subsets of a set $A$, then writing complements as $X_i' \coloneqq \!A\setminus \!X_i$, we have 
\begin{equation}
\label{complements}
    |\blacktriangle(X_1, \ldots, X_d)| = |\blacktriangle(X'_1, \ldots, X'_d)|.
\end{equation}
Note also the connection to the generalized cost function $C$ from~\eqref{cost}: we can rewrite $\blacktriangle$ as a multiset union over the distinct elements $x \in \bigcup_i X_i$, as
\[
\blacktriangle(X_1, \ldots, X_d) = \bigcup_x \: x^{\mathbf{u} \cdot (\overbrace{\scriptstyle 1, \ldots, 1}^{\mathclap{\text{$\# X_i$ containing $x$}}}, 0, \ldots, 0)},
\]
and therefore we have the cardinality
\begin{align}
    |\blacktriangle(X_1, \ldots, X_d)| &= \sum_x \mathbf{u} \cdot (\overbrace{1, \ldots, 1}^{\mathclap{\text{$\# X_i$ containing $x$}}}, 0, \ldots, 0) \nonumber \\
    &= \mathbf{u} \cdot \sum_x (\overbrace{1, \ldots, 1}^{\mathclap{\text{$\# X_i$ containing $x$}}}, 0, \ldots, 0) \label{size symm diff intermediate}\\
    &= C\big(|X_1|, \ldots, |X_d|\big), \label{size symm diff}
\end{align}
where the last equality follows from~\eqref{cost} and the fact that the sum over $x$ in~\eqref{size symm diff intermediate} equals the vector consisting of the cardinalities $|X_i|$ sorted in decreasing order.  
We arrive at the following generalization of Proposition~\ref{prop:EMDisSymmDiff}, expressing $\EMD$ in terms of Young diagrams:

\begin{proposition}
\label{prop:EMDisUSD}
For $\bla \in \H(m,n)^d$, we have \[
\EMD(\bla) = \left|\blacktriangle\big(\Y(\la^1),\ldots,\Y(\la^d)\big)\right|.
\]
\end{proposition}

\begin{proof}
Let $x^i_j$ be the $j$th data point of $\la^i$, as written out in~\eqref{histogram to data points general d}.  Recall that we constructed $\Y(\la^i)$ so that its $j$th row has length $n - x^i_j$.  Summing up the sizes of the symmetric difference over each row of the diagrams $\Y(\la^i)$, we have
\begin{align*}
    \left|\blacktriangle\big(\Y(\la^1),\ldots,\Y(\la^d)\big)\right| &= \sum_{j=1}^m \left|\blacktriangle\big([n-x^1_j], \ldots, [n-x^d_j]\big)\right|\\
    &= \sum_j \left|\blacktriangle\big([x^1_j], \ldots, [x^d_j]\big)\right|\\
    &= \sum_j C(x^1_j, \ldots, x^d_j)\\
    &= \EMD(\bla),
\end{align*}
where the last three equalities follow directly from~\eqref{complements}, \eqref{size symm diff}, and~\eqref{EMD RSK general d}, respectively.
\end{proof}

\begin{example}
\label{ex:EMD d general Young diagrams}
We revisit Example~\ref{example:RSK d general}, but this time we will calculate $\EMD$ in terms of Young diagrams, using Proposition~\ref{prop:EMDisUSD}.  Recall that $d=4$, with $n=5$ and $m=6$.  We translate each histogram into its diagram:
\begin{align*}
    \la^1 &= (4,1,1,0,0) \leadsto  111123 \leadsto 444432 \\
    \la^2 &= (3,0,0,0,3) \leadsto  111555 \leadsto 444000 \\
    \la^3 &= (0,4,2,0,0) \leadsto  222233 \leadsto 333322 \\
    \la^4 &= (1,1,2,1,1) \leadsto  123345 \leadsto 432210
\end{align*}
\[
\Y(\la^1) = \ydiagram{4,4,4,4,3,2}
 \quad 
 \Y(\la^2) = \ydiagram{4,4,4} \quad
 \Y(\la^3) = \ydiagram{3,3,3,3,2,2} 
 \quad 
 \Y(\la^4) = \ydiagram{4,3,2,2,1}
 \]
Now we superimpose these four diagrams so that they share a common upper-left box; see the two diagrams below in \eqref{diagramsexample}.  On the left-hand side of \eqref{diagramsexample}, we label each box with the number of diagrams $\Y(\la^i)$ that contain it, which produces a plane partition representing the multiset union of the four Young diagrams.  Then on the right-hand side, we represent $\blacktriangle(\Y(\la^1), \ldots, \Y(\la^4))$ by replacing each entry $k$ with $\min\{k,\:4-k\}$:
\begin{equation}\label{diagramsexample}
\ytableausetup{smalltableaux}
\begin{ytableau}
4 & 4 & 4 & 3\\
4 & 4 & 4 & 2\\
4 & 4 & 3 & 2\\
3 & 3 & 2 & 1\\
3 & 2 & 1 \\
2 & 2
\end{ytableau}
\qquad \leadsto \qquad
\begin{ytableau}
0 & 0 & 0 & *(lightgray) 1\\
0 & 0 & 0 & *(gray) 2\\
0 & 0 & *(lightgray) 1 & *(gray) 2\\
*(lightgray) 1 & *(lightgray) 1 & *(gray) 2 & *(lightgray) 1\\
*(lightgray) 1 & *(gray) 2 & *(lightgray) 1\\
*(gray) 2 & *(gray) 2
\end{ytableau}
\end{equation}
(The shading is superfluous, of course, but helps to illustrate what our generalized symmetric difference of Young diagrams looks like.)  We sum the entries of the right-hand diagram to conclude that $\EMD(\bla) = 19$,
 agreeing with the calculation in Example \ref{example:RSK d general}.
\end{example}

\subsection{Generalization of $|\D|$}

We must next generalize $|\D(\bla)|$ for $d$ histograms.  Intuitively, we want $|\D(\la^1, \ldots, \la^d)|$ to equal the minimum amount of work required to equalize all the weighted totals $t(\la^i)$.  Recall that one unit of work is equivalent to moving one data point (in any histogram) by one bin; hence one unit of work adds $\pm 1$ to the weighted total of the chosen histogram, and therefore $|\D(\bla)|$ should equal the minimum number of $\pm 1$'s required to be added to the coordinates of $\D(\bla)$ in order to equalize them.  But this is precisely the taxicab distance to the diagonal in $(\mathbb{Z}_{>0})^d$, which coincides with the generalized cost function $C$ in~\eqref{cost}.  Therefore we define
\begin{equation}
    \label{|D| is C}
    |\D(\bla)| \coloneqq C(\D(\bla))
\end{equation}
with $C$ as in~\eqref{cost}.  Note that $C$ is well-defined on $\Lambda_d$ since $C$ is invariant under addition by $(1, \ldots, 1)$. (For those readers regarding $\Lambda_d$ as the Type-A root lattice, $|\D(\bla)|$ is the distance between $\D(\bla)$ and $0$ on that lattice.)  

It will be helpful to visualize $|\D(\bla)|$ via Young diagrams, similarly to our $\EMD$ construction in~\eqref{diagramsexample}.  First, in order to use~\eqref{cost} and~\eqref{|D| is C}, we must sort the $\la^i$ in decreasing order of weighted total: that is, let $\widetilde{\bla} = (\widetilde{\la}^1, \ldots, \widetilde{\la}^d)$ be some permutation of the $\la^i$ such that $t(\widetilde{\la}^1) \geq \cdots \geq t(\widetilde{\la}^d)$.  Note that the choice of $\widetilde{\bla}$ is not unique whenever any of the $\la^i$ have the same weighted total, but clearly $|\D(\bla)|$ is independent of our choice of $\widetilde{\bla}$.  Now one by one, in order, superimpose the diagrams $\Y(\widetilde{\la}^1), \ldots, \Y(\widetilde{\la}^d)$ so that they share a common upper-left box.  The first $\lfloor d/2 \rfloor$ diagrams each contribute $1$ in every box, while the last $\lfloor d/2 \rfloor$ diagrams each contribute $-1$ in every box.  (When $d$ is odd, the median diagram contributes $0$ in every box.)

\begin{example}
\label{ex:D with Young diagrams}
    We continue Example~\ref{ex:EMD d general Young diagrams}.  Counting boxes in each of the $\Y(\la^i)$, we see that
    \[
    t(\la^1) = 21, \quad t(\la^2) = 12, \quad t(\la^3) = 16, \quad t(\la^4) = 12.
    \]
    We have two choices for $\widetilde{\bla}$ since $t(\la^2) = t(\la^4)$; we arbitrarily choose
    \[
    \widetilde{\bla} = (\la^1, \la^3, \la^2, \la^4).
    \]
    At this point, we can easily compute $|\D(\bla)| = (21 + 16) - (12 + 12) = 13$, but later we will need to understand the Young diagram construction as well, which is as follows:
    \[
    \ytableausetup{boxsize=2.3ex}
        \left(\ydiagram{4,4,4,4,3,2} + \ydiagram{3,3,3,3,2,2}\right) - \left( \ydiagram{4,4,4} + \ydiagram{4,3,2,2,1} \right)
        \quad \leadsto \quad
        \begin{ytableau}[\scriptstyle]
        0 & 0 & 0 & -1\\
        0 & 0 & 0 & 0\\
        0 & 0 & 1 & 0\\
        1 & 1 & 2 & 1\\
        1 & 2 & 1\\
        2 & 2
\end{ytableau}
    \]
    Summing up the entries of the resulting diagram, we again have $|\D(\bla)| = 13$, whereas we saw that $\EMD(\bla) = 19$.  Crucially, observe that there are exactly three boxes whose entries differ from those in the $\EMD$ diagram~\eqref{diagramsexample}, namely the top three boxes in the rightmost column: these are the boxes which belong to one of the two smaller diagrams but do not belong to one of the two larger diagrams.  This will be the idea behind Theorem~\ref{thm:EMD=D for general d}: the obstruction to equality between $\EMD$ and $|\D|$ is the failure of any of the ``positive'' diagrams to contain any of the ``negative'' diagrams. 
\end{example}
 The upshot of Examples~\ref{ex:EMD d general Young diagrams} and~\ref{ex:D with Young diagrams} is that both $|\D|$ and $\EMD$ can be computed by summing the contributions of the individual boxes in the union $U \coloneqq \bigcup_i \Y(\la^i)$, where by ``contributions'' we mean the entries of the final diagram obtained in each example.  For each box $b \in U$, let $D_b$ (resp. $E_b$) denote the contribution of $b$ to $|\D(\bla)|$ (resp., to $\EMD$).  Thus we write
 \[
 |\D(\bla)| = \sum_{b \in U} D_b, \qquad \EMD(\bla) = \sum_{b \in U} E_b.
 \]
 It is clear from the examples, and from the definitions of $C$ and of $\blacktriangle$, how to express the contributions $D_b$ and $E_b$.  First, given a sorting $\widetilde{\bla}$, we abbreviate the indicator function
 \[
 \mathbf{1}_i(b) \coloneqq \begin{cases}
     1, & b \in \Y(\widetilde{\la}^i),\\
     0, & b \notin \Y(\widetilde{\la}^i).
 \end{cases}
 \]
Then recalling the vector $\mathbf{u}$ from~\eqref{u}, we have 
\begin{equation}
\label{Db Eb}
D_b = \mathbf{u} \cdot \left(\mathbf{1}_1(b), \ldots, \mathbf{1}_d(b)\right), \qquad E_b = \mathbf{u} \cdot (\underbrace{1, \ldots, 1}_{\mathclap{\text{$\# \Y(\la^i)$ containing $b$}}}, 0, \ldots, 0).
\end{equation}
Compare the right-hand vectors in each of the two dot products above. 
 Note that they necessarily have the same number of $1$'s.  Therefore, since the coordinates of $\mathbf{u}$ are weakly decreasing, we have $E_b \geq D_b$ for all $b \in U$, and therefore $\EMD(\bla) \geq |\D(\bla)|$ in general.  We have equality if and only if, for all $b \in U$, those two right-hand vectors have the same number of $1$'s among their first $\lceil d/2 \rceil$ coordinates, since those correspond to the $1$'s in $\mathbf{u}$.

\section{Main results for $d$ histograms}
\label{section:general d results}

\subsection{Generalized relationship between $\EMD$ and $|\D|$}

Given $\bla = (\la^1, \ldots, \la^d) \in \H(m,n)^d$, let $\operatorname{m}_{\bla}$ denote the median of the weighted totals $t(\la^1), \ldots, t(\la^d)$.  Then let
\begin{align*}
    \bla^{\geq{\rm m}} & \coloneqq \{ \la^i \mid t(\la^i) \geq \operatorname{m}_{\bla}\},\\
    \bla^{\leq{\rm m}} & \coloneqq \{ \la^j \mid t(\la^j) \leq \operatorname{m}_{\bla}\}.
\end{align*}
We will speak of the ``first half'' and ``second half'' of $d$-tuples, by which we mean the first (resp., last) $\lceil d/2 \rceil$ elements.  (This includes the middle element when $d$ is odd.)  In any sorting $\widetilde{\bla}$, the histograms in the first half of $\widetilde{\bla}$   must all belong to $\bla^{\geq{\rm m}}$, while the second half  must all belong to $\bla^{\leq{\rm m}}$.
Conversely, for each $\la^i \in \bla^{\geq{\rm m}}$ and $\la^j \in \bla^{\leq{\rm m}}$, there exists a sorting $\widetilde{\bla}$ such that $\la^i$ is in the first half of $\widetilde{\bla}$, and and $\la^j$ is in the second half of $\widetilde{\bla}$.  

\begin{theorem}
\label{thm:EMD=D for general d}
    Let $\bla \in \H(m,n)^d$.  Then $\EMD(\bla) = |\D(\bla)|$ if and only if $\Y(\la^i) \supseteq \Y(\la^j)$ for each $\la^i \in \bla^{\geq {\rm m}}$ and $\la^j \in \bla^{\leq {\rm m}}$.
\end{theorem}

\begin{proof}
    First suppose that $\Y(\la^i) \supseteq \Y(\la^j)$ for each $i \in \bla^{\geq {\rm m}}$ and $j \in \bla^{\leq {\rm m}}$; we need to show that $D_b = E_b$ for all $b \in U$.  Let $\widetilde{\bla}$ be any sorting.  Let $b \in U$.  If $b$ is not contained in any of the Young diagrams of the second half of $\widetilde{\bla}$, then the second half of $(\mathbf{1}_1(b), \ldots, \mathbf{1}_d(b))$ is all 0's, and hence $D_b = E_b$.  Hence suppose that $b$ is contained in some diagram $\Y(\widetilde{\la}^{k})$ where $\widetilde{\la}^{k}$ is in the last half of $\widetilde{\bla}$.  Then $\widetilde{\la}^{k}\in \bla^{\leq{\rm m}}$, and thus $b$ is contained in all of the Young diagrams of the first half of $\widetilde{\bla}$.  Therefore the first half of $(\mathbf{1}_1(b), \ldots, \mathbf{1}_d(b))$ is all $1$'s, and hence $D_b = E_b$.

    To prove the converse, suppose that there is some $\la^i \in \bla^{\geq{\rm m}}$ and $\la^j \in \bla^{\leq{\rm m}}$ such that $\Y(\la^i) \not \supseteq \Y(\la^j)$.  Then there exists $b \in U$ such that $b \in \Y(\la^j)$ but $b \notin \Y(\la^i)$.  Choose a sorting $\widetilde{\bla}$ in which $\la^i$ is in the first half and $\la^j$ is in the second half.  Then $(\mathbf{1}_1(b), \ldots, \mathbf{1}_d(b))$ has a 0 in the first half and a 1 in the second half, and hence $E_b > D_b$.
\end{proof}

It seems that it will be difficult to enumerate the plane partitions corresponding to the families of diagrams $\Y(\la^i)$ described in Theorem~\ref{thm:EMD=D for general d}.  For now, then, we leave the full generalization of Theorem~\ref{theorem:EMDequalsD} as an open problem:
\begin{problem}
\label{prob:EMD=D general d}
    Determine an explicit formula for $\mathbb{P}(\EMD(\bla) = |\D(\bla)|)$ that holds for arbitrary $d$.
\end{problem}
 In the case $d=3$, however, we observe that Theorem~\ref{thm:EMD=D for general d} forces the following: $\EMD(\bla) = |\D(\bla)|$ if and only if the three diagrams $\Y(\la^i)$ form a chain with respect to containment.  Therefore, as a first step toward completing Problem~\ref{prob:EMD=D general d}, we solve the $d=3$ case below.

\begin{proposition}
\label{prop:EMD=D for d=3}
    Let $\bla \in \H(m,n)^3$ be chosen uniformly at random. 
 Then the probability that $\EMD(\bla) = |\D(\bla)|$ is
 \begin{equation}
    \label{formula EMD=D for d=3}
    \frac{12 (m+n)^2 (m+n+1)}{n^2 (n+1) (m+1)^2 (m+2)} - \frac{6(m+n)m!}{n(m+1)(n)_m} + \left(\frac{m!}{(n)_m}\right)^2.
\end{equation}
\end{proposition}

\begin{proof}
    By Theorem~\ref{thm:EMD=D for general d}, the desired probability equals the number of ordered triples in $\scrY(m, n-1)^3$ which form a chain with respect to containment, divided by $\# \scrY(m, n-1)^3$.   We claim that the number of ordered triples in $\scrY(m, n-1)$ forming a chain equals
    \begin{equation}
        \label{3-chains}
        6\PP(m,\:n-1,\:3) - 6\PP(m,\:n-1,\:2) + \PP(m,\:n-1,\:1).    
    \end{equation}
This can be seen as follows.  By~\eqref{bijection PP's}, the number of \emph{unordered} triples $Y_1 \supseteq Y_2 \supseteq Y_3$ equals $\PP(m,\:n-1,\:3)$.  To obtain the number of such \emph{ordered} triples, we multiply by $3! = 6$, but this overcounts: the cases in which $Y_1 = Y_2$, and the cases in which $Y_2 = Y_3$, should have been multiplied only by $3$, while the cases in which $Y_1 = Y_2 = Y_3$ should have been multiplied only by 1.  The cases in which $Y_1 = Y_2$ correspond to the plane partitions without 1's; the cases in which $Y_2 = Y_3$ correspond to the plane partitions without 2's.  Both of these cases are equinumerous with the plane partitions without 3's, of which there are $\PP(m,\:n-1,\:2)$.  Hence we subtract $3$ copies of $\PP(m,\: n-1, \: 2)$ for each of the two cases ($Y_1 = Y_2$ and $Y_2 = Y_3$).  But doing this also subtracts all $3 \cdot 2 = 6$ copies of the plane partitions containing only 3's, which correspond to the case $Y_1 = Y_2 = Y_3$; hence we need to add them back, and there are $\PP(m,\:n-1,\:1)$ of these.  We thus obtain~\eqref{3-chains}, as claimed.

Now substituting the specializations~\eqref{PP1}--\eqref{PP3} into the expression~\eqref{3-chains}, and then dividing by the third power of $\# \scrY(m, n-1) = \PP(m,\: n-1, \: 1) = \frac{(n)_m}{m!}$, we find that our desired probability equals
\[
\left(\frac{12(n)_m (n+1)_m (n+2)_m}{m! (m+1)! (m+2)!} - \frac{6(n)_m (n+1)_m}{m!(m+1)!} + \frac{(n)_m}{m!}\right) \Bigg/ \left(\frac{(n)_m}{m!}\right)^3
\]
which we simplify as
\[
\frac{12 (m!)^2 (n+1)_m (n+2)_m}{(m+1)! (m+2)! (n)_m^2} - \frac{6(m!)^2 (n+1)_m}{(m+1)!(n)_m^2} + \frac{(m!)^2}{(n)_m^2}
\]
which we further simplify in the form~\eqref{formula EMD=D for d=3}.   
\end{proof}

Recall the continuous analogues $\widehat{\EMD}$ and $\widehat{\D}$ that we defined before Corollary~\ref{cor:continuous}.  Just as before, $\widehat{\EMD}$ and $\widehat{\D}$ have natural analogues on $\mathcal{P}(n)^d$, that is, on $d$-tuples of probability distributions on the set $[n]$.

\begin{corollary}
    Let $\mathbf{X} \in \mathcal{P}(n)^3$ be chosen uniformly at random.  Then
    \[
    \mathbb{P}\Big( \widehat{\EMD}(\mathbf{X}) = |\widehat{\D}(\mathbf{X})| \Big) = \frac{12}{n^2(n+1)}.
    \]
\end{corollary}

\begin{proof}
    We fix $n \geq 2$ and take the limit of~\eqref{formula EMD=D for d=3} as $m \rightarrow \infty$.  Note that 
    \[
    \lim_{m \rightarrow \infty} \frac{m!}{(n)_m} \leq \lim_{m\rightarrow \infty} \frac{m!}{(2)_m} = \lim_{m \rightarrow \infty} \frac{1}{m+1} = 0,
    \]
    and therefore the second and third terms in~\eqref{formula EMD=D for d=3} vanish as $m \rightarrow \infty$.  The limit of the first term is $12/n^2(n+1)$.
\end{proof}

\subsection{Expected value of $|\D|$}

The following lemma generalizes the Laurent polynomial $f(q)$ from Theorem~\ref{theorem:expd}:

\begin{lemma}
\label{lemma:D distribution}
The generating function for $\D$ on $\H(m,n)^d$ is given by 
\begin{equation}
\label{D distribution}
   f(q_1, \ldots, q_{d-1}) = \left( \prod_{i=1}^{d-1} \genfrac[]{0pt}{}{m+n-1}{m}_{q_{i}} \right) \genfrac[]{0pt}{}{m+n-1}{m}_{(q_1 \cdots q_{d-1})^{-1}}.
\end{equation}
Specifically, the coefficient of $q_1^{t_1} \cdots q_{d-1}^{t_{d-1}}$ in $f$ equals the number of elements $\bla \in \H(m,n)^d$ such that $\D(\bla) = [t_1,\ldots,t_{d-1},0]$.
\end{lemma}

\begin{proof}
By \eqref{qbinyoungfact}, the coefficient of $q_1^{t_1}\cdots q_d^{t_d}$ in the expansion of
\[
\prod_{i=1}^{d} \genfrac[]{0pt}{}{m+n-1}{m}_{q_{i}}
\]
equals the number of elements $\bla = (\la^1,\ldots,\la^d) \in \H(m,n)^d$ such that $t(\la^i) = t_i$ for each $i=1,\ldots,d$.  Carrying out the substitution $q_d = (q_1 \cdots q_{d-1})^{-1}$ then identifies terms whose ratio is some power of $q_1 \cdots q_d$, and leaves the expansion expressed in terms of $q_1,\ldots,q_{d-1}$ only.  Hence the coefficient of $q_1^{t_1} \cdots q_{d-1}^{t_{d-1}}$ now equals the number of $\bla$'s such that $\D(\bla) = [t_1,\ldots,t_{d-1},0]$.
\end{proof}

We write $q^{\mathbf t} \coloneqq q_1^{t_1} \cdots q_{d-1}^{t_{d-1}}$.  Given an exponent vector $\mathbf t = (t_1, \ldots, t_{d-1})$, let $z(\mathbf t)$ be the vector that tallies the multiplicities of the distinct numbers in the list $(t_1, \ldots, t_{d-1}, 0)$, in descending order.  (Note the trailing zero.)  For example, if $\mathbf t = (6,4,4,3,4,0,4,4,3)$, then $z(\mathbf t) = (5,2,2,1)$, since there are $5$ fours, $2$ threes, $2$ zeros, and $1$ six.  Note that the components in $z(\mathbf t)$ necessarily sum to $d$.  We will write multinomial coefficients in the form $\binom{d}{z(\mathbf t)}$.  In the previous example, this yields
\[
\binom{d}{z(\mathbf t)} = \binom{10}{5,2,2,1} = \frac{10!}{5! 2! 2! 1!} = 7560.
\]
Recall that this multinomial coefficient is the number of distinguishable ways to permute the coordinates $(t_1, \ldots, t_{d-1}, 0)$.  We use this to generalize the plussing operator $[ \quad ]_+$ from Theorem~\ref{theorem:expd}, by defining
\begin{equation}
    \label{generalized plussing}
[q^{\mathbf t}]_+ \coloneqq 
\begin{cases}
    \binom{d}{z(\mathbf t)} q^{\mathbf t}, & t_1 \geq \cdots \geq t_{d-1} \geq 0,\\
    0 & \text{otherwise}
\end{cases}
\end{equation}
on monomials, and extending by linearity.  

\begin{theorem}
\label{thm:expval |D| general d}
    Let $\bla \in \H(m,n)^d$ be chosen uniformly at random.  Then
    \[    \mathbb{E}\big(|\D(\bla)|\big) = \frac{g'(1)}{\binom{m+n-1}{m}_,^d}
    \]
    with $g(q) \coloneqq [f]_+ \Big|_{q_i = q^{u_i}}$, and $\mathbf{u} = (u_1, \ldots, u_d)$ as in~\eqref{u}.
\end{theorem}

\begin{proof}
    In this proof, we will identify exponent vectors with elements of $\Lambda_d$, via $\mathbf t \leftrightarrow [t_1, \ldots, t_{d-1},0]$.  If $t_1 \geq \cdots \geq t_{d-1}$, then it follows from~\eqref{cost} and~\eqref{|D| is C} that applying the substitutions $q_i = q^{u_i}$ to a monomial results in the following:
    \[
    q^{\mathbf t} = q_1^{t_1} \cdots q_{d-1}^{t_{d-1}} \leadsto q^{\mathbf{u} \cdot (t_1, \ldots, t_{d-1},0)} = q^{C(t_1, \ldots, t_{d-1},0)} = q^{|\mathbf t|}.
    \]
    Now consider the symmetric group on $d$ letters, acting on $\Lambda_d$ by permutation of coordinates.  By \eqref{cost} and \eqref{|D| is C}, we observe that $|\mathbf{t}|$ is constant on the orbit of $\mathbf{t}$.  We have therefore defined the generalized plussing operator in~\eqref{generalized plussing} to pick out a unique representative from the orbit of each $\mathbf{t} \in \{\D(\bla) \mid \bla \in \H(m,n)^d\}$, and we have scaled by the multinomial coefficient, which is the size of that orbit.  Hence by Lemma~\ref{lemma:D distribution} we have
    \begin{align*}
    g(q) &= \sum_{t_1 \geq \cdots \geq t_{d-1} \geq 0} \#\{\bla \mid \D(\bla) = \mathbf{t}\} \cdot \#(\text{orbit of $\mathbf t$}) \:q^{|\mathbf{t}|}\\
    &= \sum_{\mathbf{t} \in \{\D(\bla) \mid \bla \in \H(m,n)^d\}} \#\{\bla \mid \D(\bla) = \mathbf{t}\} \: q^{|\mathbf t|}\\
    &= \sum_{\bla \in \H(m,n)^d} q^{|\D(\bla)|}.
    \end{align*}
Thus $g'(1) = \sum_{\bla} |\D(\bla)|$, and the result follows upon division by $\#\H(m,n)^d$.
\end{proof}

By means of Theorem~\ref{thm:expval |D| general d}, we are able to use software to compute the expected value of $|\D|$ quite efficiently; for example, Mathematica takes about 10 seconds to determine that the expected value is $30.0176$ on $\H(20,10)^3$, which contains just over $10^{21}$ triples.

\begin{figure}[h]
    \centering
\includegraphics[width=0.65\textwidth]{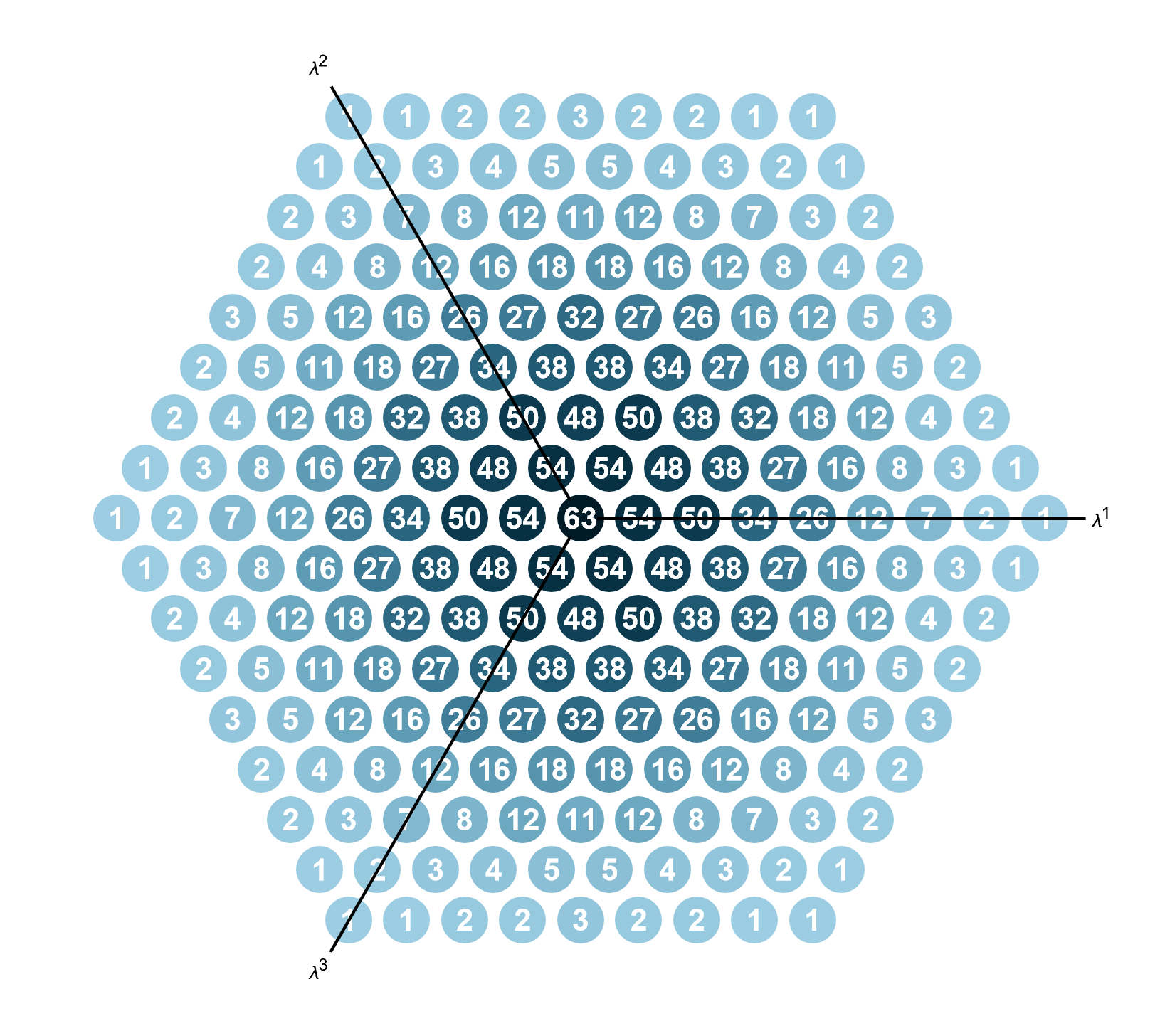}
    \caption{Distribution of $\D$-values on $\H(4,3)^3$, accompanying Example~\ref{fig:example d3}.  The hexagonal lattice is the image of $\D$ inside $\Lambda_3$.}
    \label{fig:example d3}
\end{figure}

\begin{example}
    It is instructive to visualize the distribution of $\D$-values when $d=3$, since this can be done on a two-dimensional plot; see Figure~\ref{fig:example d3}.  We take $m=4$ and $n=3$.  The hexagonal lattice is the image of $\D$ inside $\Lambda_3$. 
    To find an element $[t_1, t_2, t_3]$, begin at the origin $[0,0,0]$; move $t_1$ units along the $\la^1$-axis, then $t_2$ units along the $\la^2$-axis, then $t_3$ units along the $\la^3$-axis. 
    The set $\H(4,3)^3$ contains $3375$ histogram triples, and these triples are distributed on the lattice according to their $\D$-value. Hence the figure is a visualization of the $\D$-generating function $f$ in Lemma~\ref{lemma:D distribution}: for example, the coefficient of $q_1^3 q_2^6$ is $11$, which we confirm by finding the ``$11$'' labeling the element $[3,6,0]$ on the lattice.  The $|\D|$-generating function $g(q)$ in Theorem~\ref{thm:expval |D| general d} is also easy to visualize, since $|\D|$ is just the distance from the origin on the lattice.  Therefore $|\D|$ increases by one on each successive concentric hexagon.  We verify this by computing $g(q)$ as outlined above:
    \[
    g(q) = 63 + 324 q + 588 q^2 + 660 q^3 + 672 q^4 + 480 q^5 + 348 q^6 + 156 q^7 + 84 q^8.
    \]
    The reader can check that the coefficients are indeed the sums along each concentric hexagon in Figure~\ref{fig:example d3}.  (See Section~\ref{section:RepThy} for an interpretation as a weight diagram for the Lie algebra $\fraksl(3,\mathbb{C})$.)
\end{example}

\begin{remark}
    The plussing operator which we used in Theorem~\ref{theorem:expd} is not quite the $d=2$ specialization of the version defined in~\eqref{generalized plussing}; if it were, it would send each nonconstant monomial $q^t$ to $\binom{2}{1,1} q^t = 2q^t$, while sending the constant $1$ to $\binom{2}{2} = 1$.  Although adopting the multinomial convention  from~\eqref{generalized plussing} would indeed have streamlined the proof of Theorem~\ref{theorem:expd}, we chose in the $d=2$ exposition to retain the standard plussing operator that appears elsewhere in the literature.  This explains the extra factor of $2$ appearing in the statement of Theorem~\ref{theorem:expd}, when compared to Theorem~\ref{thm:expval |D| general d}.
\end{remark}

\section{The view from representation theory}
\label{section:RepThy}

We conclude by pointing out the representation-theoretic perspective on our results above.  (For readers unfamiliar with the subject, a glance at the diagrams in \cite[Lectures 12--13]{fh} will suffice to see the connections we describe below.) Recall that the generalized weighted difference $\D$ has codomain $\Lambda_d \coloneqq \mathbb{Z}^d / \mathbb{Z}(1, \ldots, 1)$.  This $\Lambda_d$ has the structure of both the root lattice and the integral weight lattice of the Lie algebra $\fraksl_d = \fraksl(d,\mathbb{C})$, i.e., of the root system ${\rm A}_{d-1}$.  

On one hand, we can identify $\Lambda_d$ with the integral weight lattice by using the standard ``epsilon'' coordinates, where as usual $\varepsilon_i$ is the linear functional sending a diagonal matrix in the Cartan subalgebra of $\fraksl_d$ to its $i$th diagonal entry: in this way, $[t_1, \ldots, t_d] = t_1 \varepsilon_1 + \cdots + t_d \varepsilon_d$.  From this perspective, $\D$ sends a $d$-tuple of histograms to an integral weight.  In particular, our generalized plussing operator \eqref{generalized plussing} retains only the \emph{dominant} weights appearing in the expansion of $f$ in Lemma~\ref{lemma:D distribution}, and the multinomial coefficient $\binom{d}{z(\mathbf{t})}$ is the size of the orbit of the weight $\mathbf{t}$ under the action of the Weyl group.  

On the other hand, an interesting question arises if we instead view $\Lambda_d$ as the root lattice. 
 For $i = 1, \ldots, d-1$, define the simple root $\alpha_i \coloneqq \varepsilon_i - \varepsilon_{i+1}$ as usual, and set $\alpha_0 \coloneqq \varepsilon_d - \varepsilon_1$, which is sometimes called the \emph{affine root}.  We then identify $\Lambda_d$ with the root lattice via the identification $[t_1, \ldots, t_d] = t_1 \alpha_1 + \cdots + t_{d-1} \alpha_{d-1} + t_d \alpha_0$. Now in Lemma~\ref{lemma:D distribution}, by setting $q_i = e^{\alpha_i}$, it turns out that $f$ is a virtual character of $\fraksl_d$.  (See details in~\cite[Thm.~6.14]{EricksonThesis}.)  Thus, the distribution of $\D$-values on $\H(m,n)^d$ can be viewed as the weight diagram of this virtual representation (see Figure~\ref{fig:example d3}).  This raises an intriguing question:
 
\begin{problem}\label{problem:decompose} How does the $\D$-distribution on $\H(m,n)^d$ decompose into irreducible $\fraksl_d$-characters?
\end{problem}

\begin{figure}[htbp]
\captionsetup[subfigure]{labelformat=empty}
 
\begin{subfigure}[t]{0.32\textwidth}
\includegraphics[width=\linewidth]{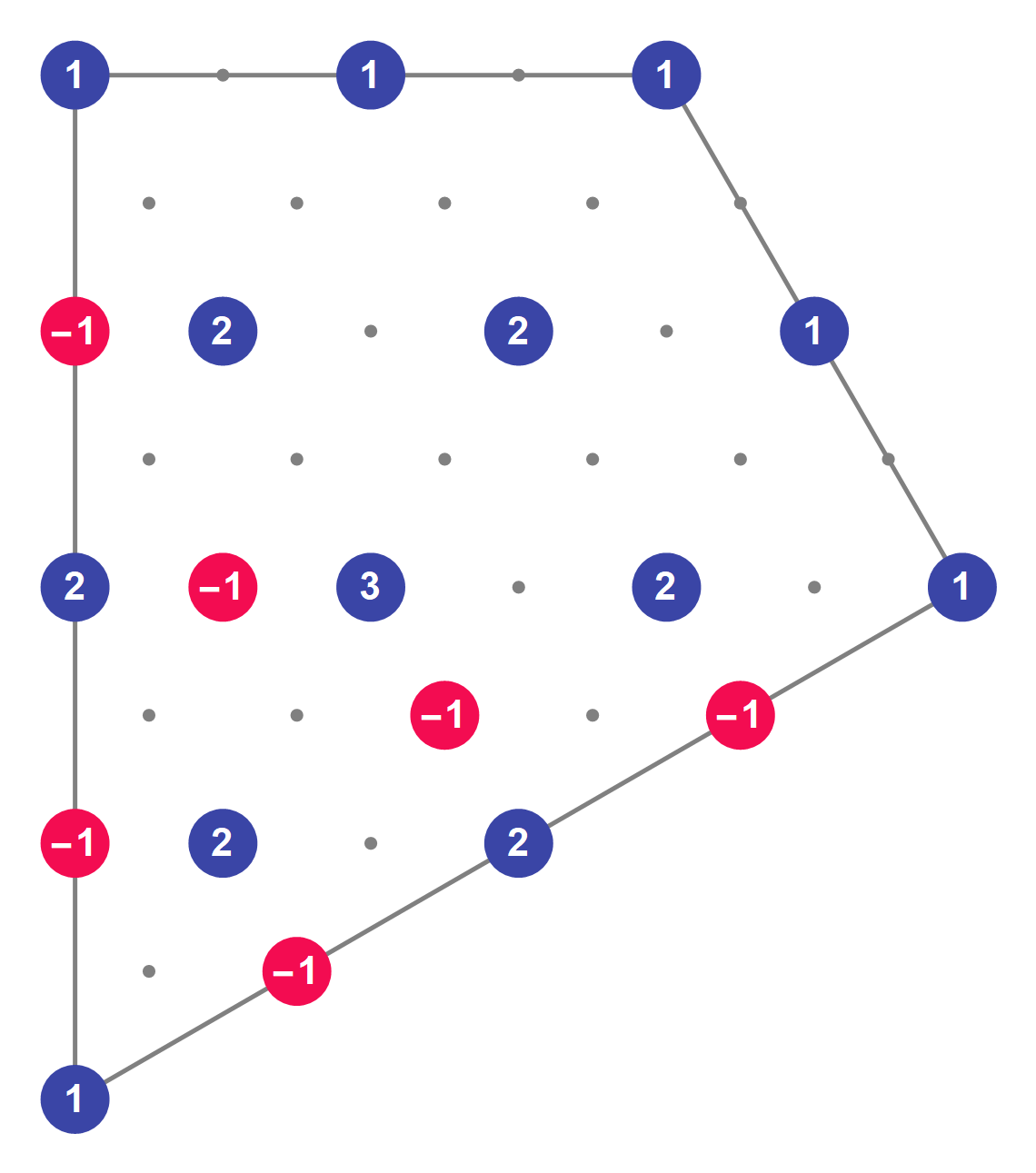} 
\caption{$m=2$, $n=5$}
\end{subfigure}
\hfill
\begin{subfigure}[t]{0.32\textwidth}
\includegraphics[width=\linewidth]{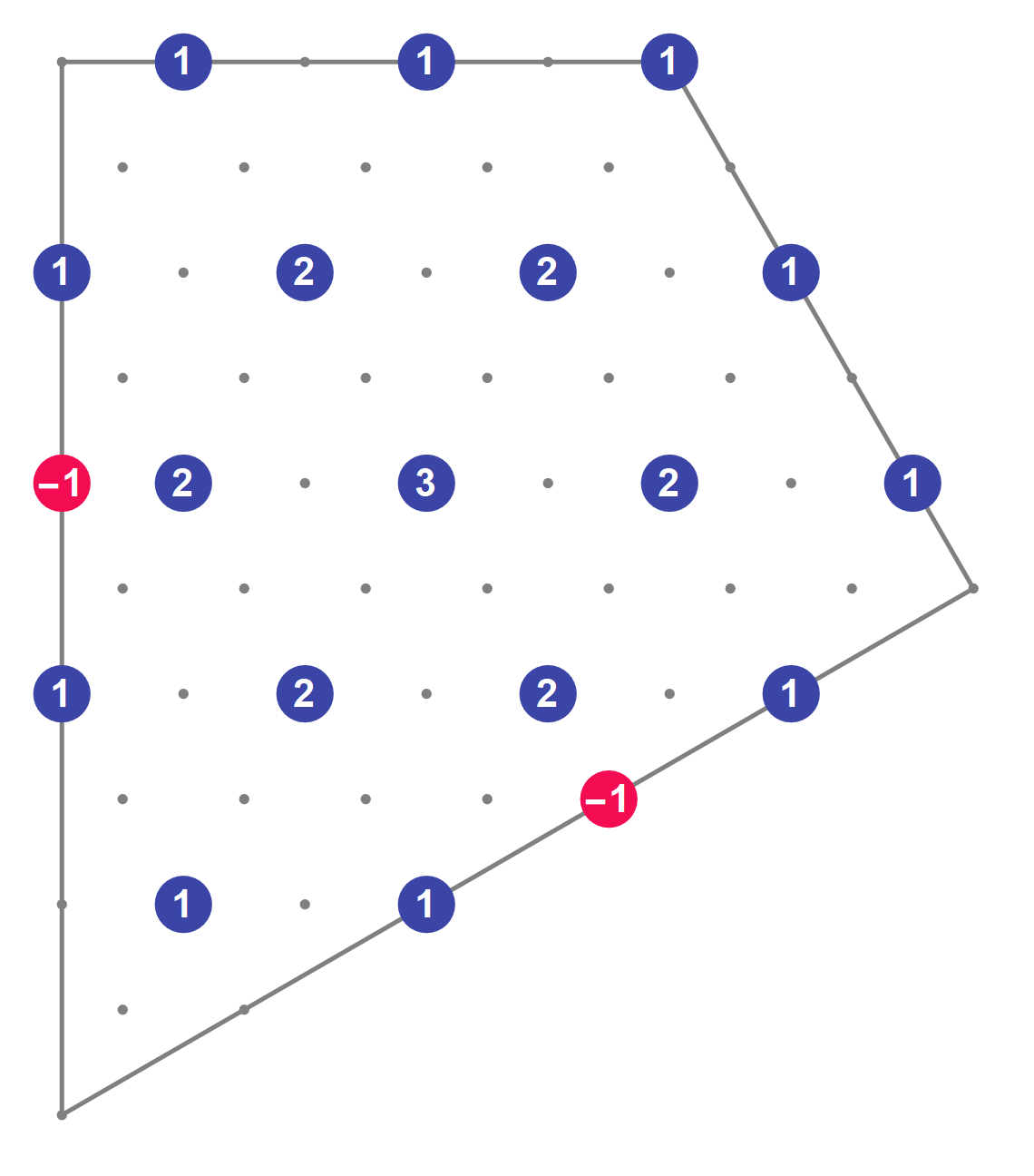} 
\caption{$m=2$, $n=6$}
\end{subfigure}
\hfill
\begin{subfigure}[t]{0.32\textwidth}
\includegraphics[width=\linewidth]{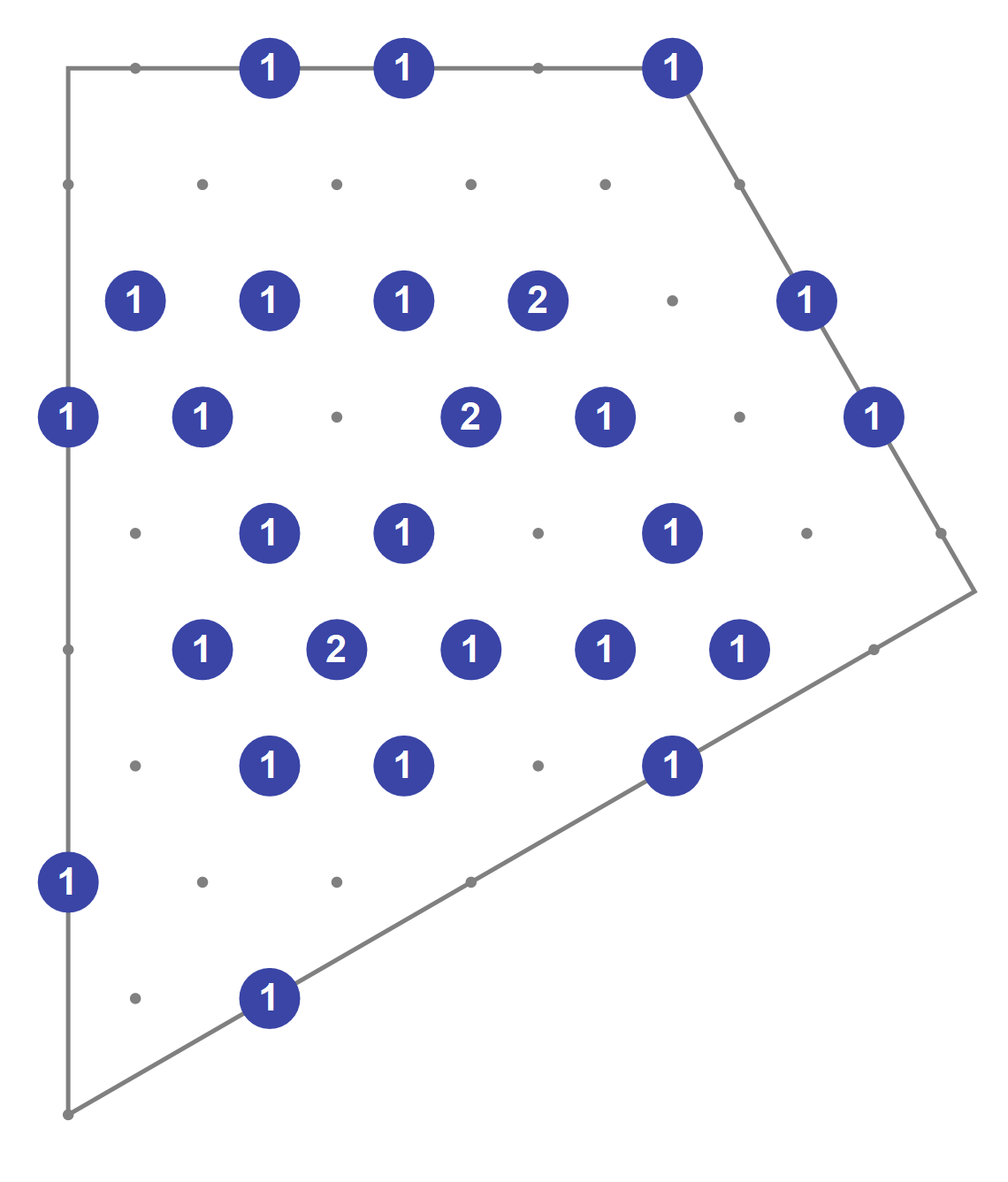} 
\caption{$m=3$, $n=4$}
\end{subfigure}

\vspace{1cm}

\begin{subfigure}[t]{0.32\textwidth}
\includegraphics[width=\linewidth]{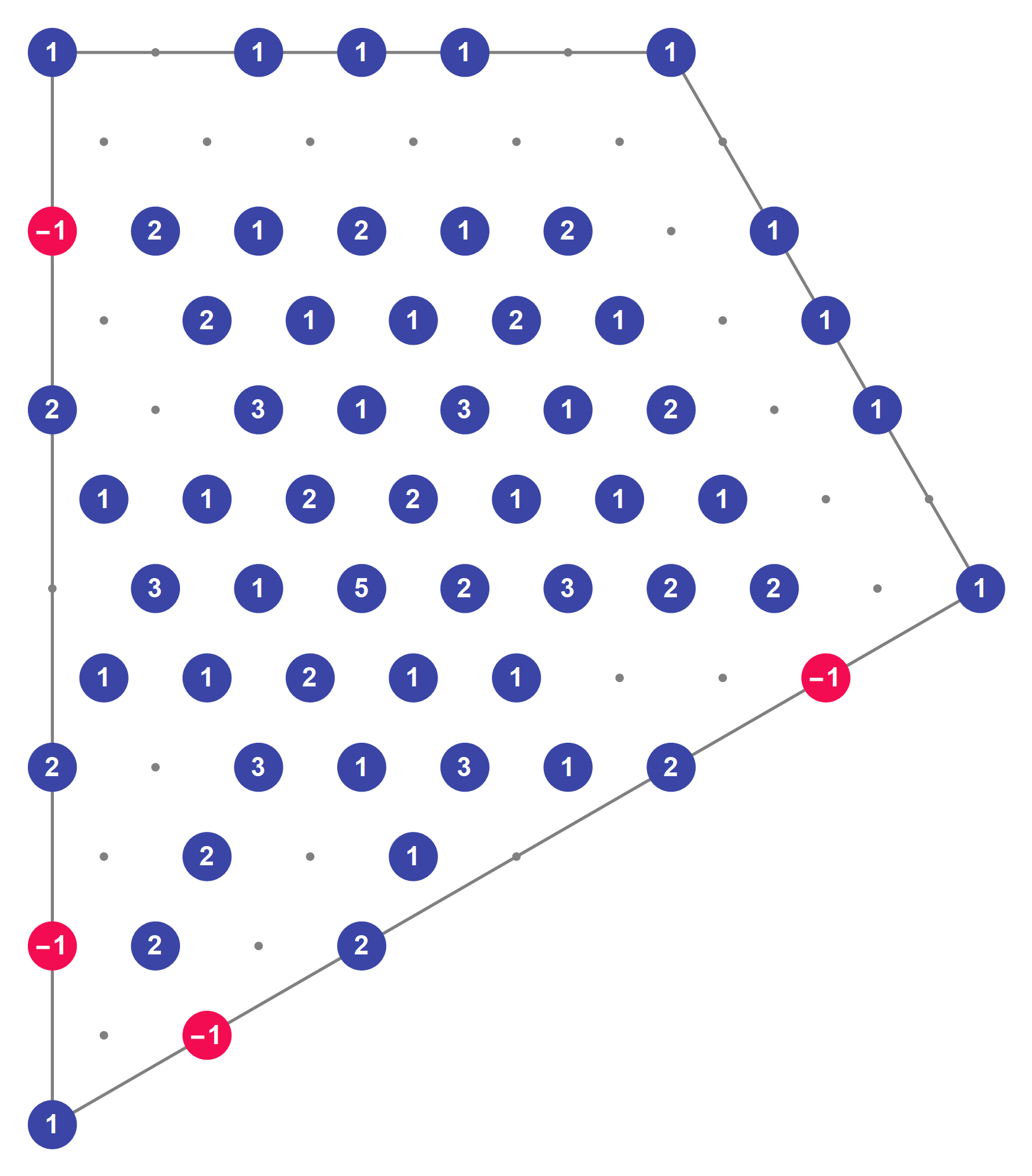} 
\caption{$m=3$, $n=5$}
\end{subfigure}
\hfill
\begin{subfigure}[t]{0.32\textwidth}
\includegraphics[width=\linewidth]{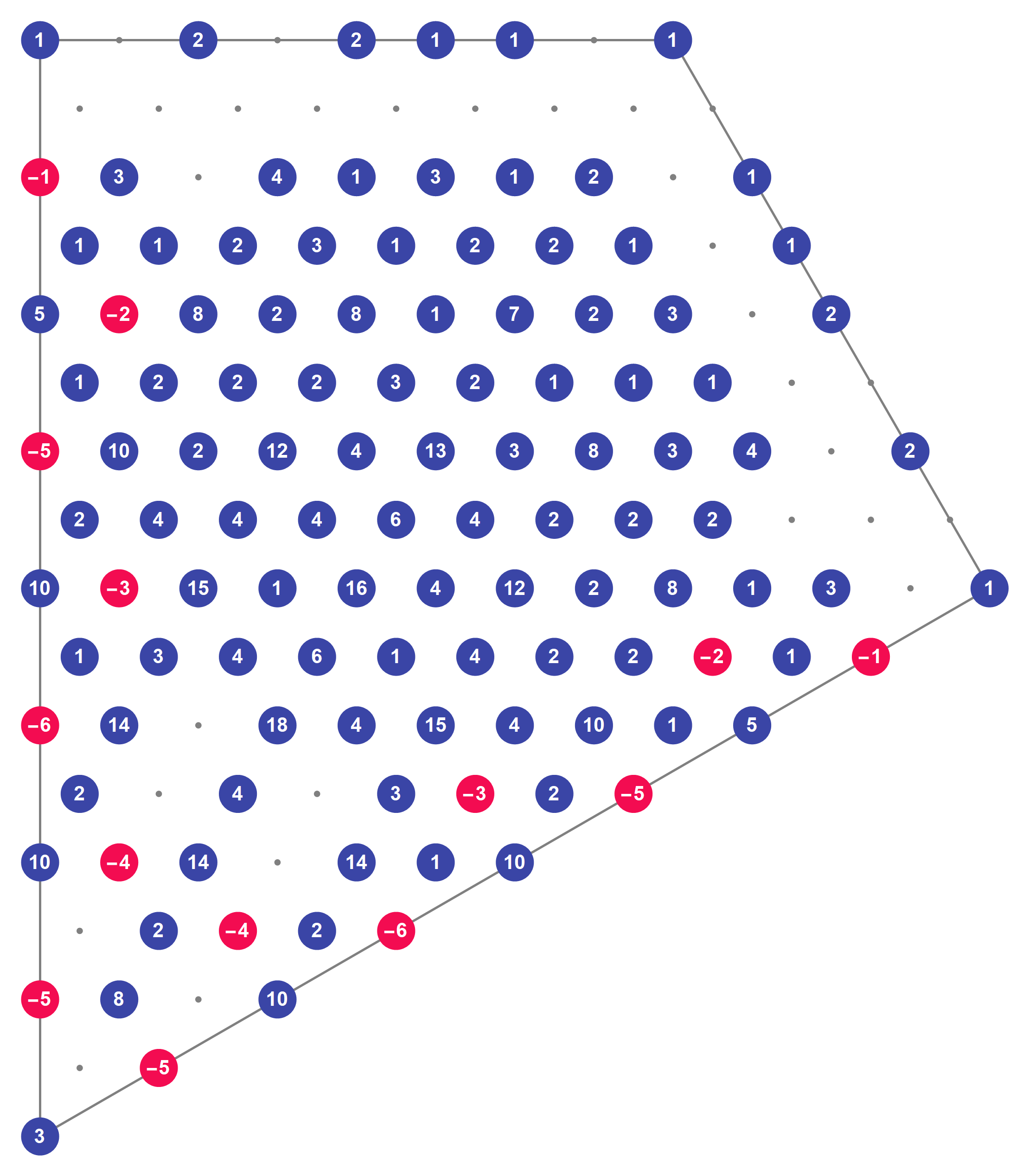} 
\caption{$m=4$, $n=5$}
\end{subfigure}
\begin{subfigure}[t]{0.32\textwidth}
\includegraphics[width=\linewidth]{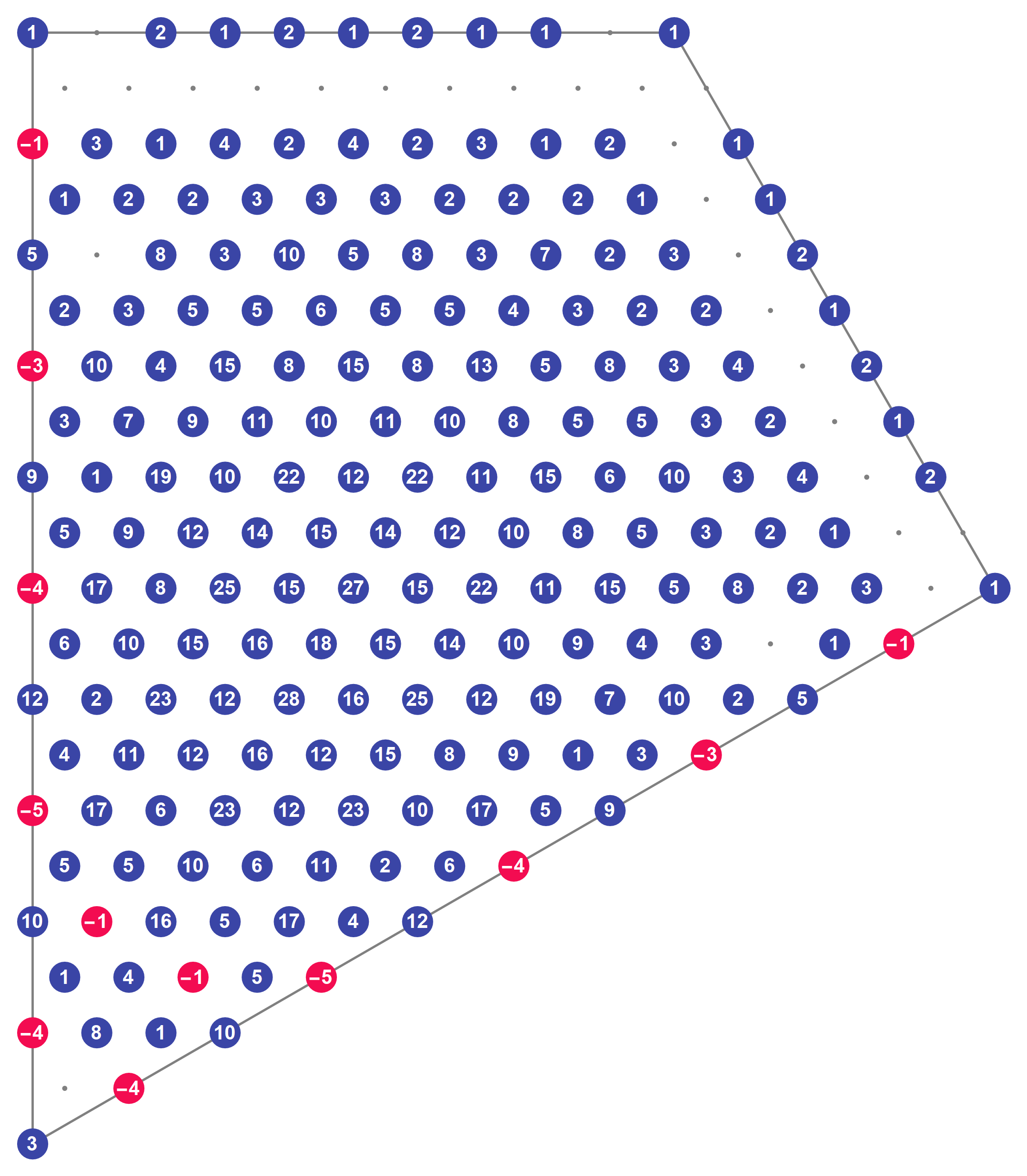} 
\caption{$m=4$, $n=6$}
\end{subfigure}

\caption{An illustration of Problem~\ref{problem:decompose}, namely, the decomposition of $f$ (Lemma~\ref{lemma:D distribution}) into irreducible $\fraksl_d$-characters.  Here we take $d=3$, and in each plot we show the dominant chamber of the root lattice for $\fraksl_3$.  Each irreducible character is plotted (with multiplicity) at its highest weight.  Since $f$ is a virtual character, positive multiplicities are colored in blue and negative multiplicities in red.}
\label{fig:d3irreps}
\end{figure}

We have carried out the decomposition in Mathematica for $d=3$, where the weight diagram is two-dimensional, for the first few values of $m$ and $n$; see Figure~\ref{fig:d3irreps}. We programmed these decompositions by multiplying the character $f$ by the Weyl denominator $(1-q_1^{-1})(1-q_2^{-1})(1-(q_1q_2)^{-1})$, and then plotting the dominant weights in the resulting expansion. Since finite-dimensional irreducible representations are indexed by their highest weight, the plots in Figure~\ref{fig:d3irreps} convey the same information as do the $\D$-distributions such as Figure~\ref{fig:example d3}, but with a much smaller and sparser diagram.  As for Problem~\ref{problem:decompose}, certain patterns seem to be emerging from these first few diagrams for $d=3$; we would like in the future to be able to write down a decomposition for all $d$, $m$, and $n$.

\bibliographystyle{alpha}
\bibliography{cas-refs}

\end{document}